\documentclass[10pt]{conm-p-l}

\usepackage{amsmath}
\usepackage{amsthm}
\usepackage{amssymb}
\usepackage{amscd}
\usepackage{diagrams}

\theoremstyle{plain}

\newtheorem{thm}[subsection]{Theorem}

\newtheorem{cor}[subsection]{Corollary}

\theoremstyle{remark}
\newtheorem{rmk}[subsection]{Remark}

\def\EE{\EE}

\def\TT{\mathbb{T}}
\def\NN{\mathbb{N}}
\def\CC{\mathbb{C}}

\def\EE{\mathcal{E}}
\def\KK{\mathcal{K}}

\def\Cat{\mathrm{Cat}}

\def\AA{\mathbb{A}}
\def\BB{\mathbb{B}}

\def\SS{\mathcal{S}}

\def\WW{\mathrm{W}}
\def\OO{\mathrm{O}}

\def\Coll{\mathrm{Coll}}
\def\Oper{\mathrm{Oper}}

\def\Sets{\mathrm{Sets}}

\def\Ass{\mathcal{A}ss}
\def\op{\mathrm{op}}
\def\Top{\mathrm{Top}}
\def\Cat{\mathrm{Cat}}
\def\Hom{\mathrm{Hom}}

\def\Alg{\mathrm{Alg}}
\def\End{\mathrm{End}}

\def\Ho{\mathbf{Ho}}
\def\Sg{\Sigma}
\def\sg{\sigma}
\def\eps{\epsilon}
\def\al{\alpha}
\def\be{\beta}
\def\lda{\lambda}
\def\ito{\rightarrowtail}

\def\lra{\longrightarrow}

\def\lrto{\leftrightarrows}
\def\rlto{\rightleftarrows}
\def\dto{\rightrightarrows}

\def\eqv{\overset{\sim}{\longrightarrow}}

\def\cli{\varinjlim}

\def\inp{in}
\def\Aut{\mathrm{Aut}}
\def\Ku{\underline{\KK}}
\def\Diag{\mathrm{Diag}}
\def\CohDiag{\mathrm{CohDiag}}
\def\Gr{\mathrm{Gr}}
\def\PLie{\mathrm{PLie}}
\def\oo{\circ}

\title[Resolution of coloured operads]{Resolution of coloured operads and rectification of homotopy algebras}

\author{Clemens Berger and Ieke Moerdijk}

\date{1 December 2005. Revised 22 May 2006.}


\begin{document}
\maketitle
\begin{center}\emph{Dedicated to Ross Street on the occasion of his 60th birthday}\end{center}

\begin{abstract}We provide general conditions under which the algebras for a coloured operad in a monoidal model category carry a Quillen model structure, and prove a Comparison Theorem to the effect that a weak equivalence between suitable such operads induces a Quillen equivalence between their categories of algebras. We construct an explicit Boardman-Vogt style cofibrant resolution for coloured operads, thereby giving a uniform approach to algebraic structures up to homotopy over coloured operads. The Comparison Theorem implies that such structures can be rectified.\end{abstract}\vspace{2ex}

Many algebraic stuctures are parametrised by operads, and deformations of such structures are then controlled by suitable resolutions of operads. Early examples in topology are the $\WW$-construction of Boardman and Vogt \cite{BV}, and the explicit $A_\infty$- and $E_\infty$- operads of Stasheff \cite{Sta}, May \cite{M} and Boardman-Vogt \cite{BV}. For operads in chain complexes, one often uses cobar-bar resolutions (cf. Ginzburg-Kapranov \cite{GK}, Kontsevich-Soibelman \cite{KS}), or the smaller Koszul resolutions (cf. Fresse \cite{F}). A coherent theory of such resolutions is provided by a Quillen homotopy theory of operads, as established by Hinich \cite{H} for operads in chain complexes, by Rezk \cite{R} for simplicial operads, and by Spitzweck \cite{Sp} and the authors \cite{BM1} for operads in general monoidal model categories. When such a homotopy theory is available, one can define for any operad $P$ the notion of $P$-algebra up to homotopy - or \emph{homotopy $P$-algebra} - in a homotopy invariant way as an algebra over a cofibrant resolution of $P$. Indeed, all important instances of homotopy $P$-algebras occurring in the literature are of this form. For many operads $P$, it can moreover be proved that homotopy $P$-algebras, defined in this way, can be \emph{rectified}, in the sense of being weakly homotopy equivalent to strict $P$-algebras. One of the first instances of this phenomenon is the well known fact that in topology, any $A_\infty$-space is weakly homotopy equivalent to a topological monoid (Stasheff \cite{Sta}, Boardman-Vogt \cite{BV}). The model theoretic framework of \cite{BM1} provides a general rectification result which includes many of the known cases.

In categories where the tensor product is the cartesian product, there is another approach to algebraic structures, based on Lawvere's notion of algebraic theory. For such an algebraic theory $T$, this gives rise to a somewhat different notion of homotopy $T$-algebra (\emph{not} based on a resolution of $T$), for which a rectification result has been proved in the context of simplicial sets by  Badzioch \cite{Bad} and Bergner \cite{Ber}. 

One purpose in this paper is to prove a multi-sorted version of such a rectification theorem, based on the notion of coloured operad. This concept goes back to Boardman and Vogt \cite{BV}, and also occurs in homotopy theory under the name \emph{multicategory} (cf. e.g. Elmendorf-Mandell \cite{EM}). The precise definition and some typical examples of coloured operads and their algebras will be reviewed in Section 1 of this paper. These examples include bimodules over monoids, enriched (e.g. topological, simplicial, or differential graded) categories, diagrams on a fixed such enriched category, morphisms between algebras over a given operad $P$, and many more. In particular, operads themselves are also algebras for a suitable coloured operad. We will discuss how the homotopy theory of operads developed in \cite{BM1,BM2} extends to coloured operads. The results will be formulated and proved in a general monoidal model category, so as to provide a uniform treatment of operads and their algebras in a variety of contexts, such as spaces, simplicial sets, chain complexes, spectra, and so on. We will prove that under certain conditions, the category of $P$-algebras for a coloured operad $P$ in a monoidal model category carries a Quillen model category structure (cf. Theorem \ref{main0}). In Section 4, we will prove a general Comparison Theorem (Theorem \ref{comparison}), which provides sufficient conditions for a weak equivalence between coloured operads to induce a Quillen equivalence between the corresponding categories of algebras. This Comparison Theorem plays a central role in the applications.

Important instances of such weak equivalences between operads are provided by coloured versions of various types of resolutions mentioned above. We will present one such resolution in detail; it is an extension of the Boardman-Vogt construction \cite{BV}. Our construction applies to an arbitrary coloured operad $P$ in a monoidal model category $\EE$ possessing a suitable interval $H$, and provides a functorial cofibrant resolution $\WW(H,P)\to P$ under some mild conditions on $P$ (cf. Theorem \ref{W}). The notion of homotopy $P$-algebra is then captured by that of a $\WW(H,P)$-algebra. For example, we obtain in this way explicit definitions of notions like ``module up to homotopy over an $A_\infty$-algebra'', ``operad up to homotopy'', etc. If the conditions of the Comparison Theorem are satisfied, we deduce as a corollary that there is a Quillen equivalence \begin{gather}\label{rectif}(W(H,P)\mathrm{-algebras})\sim(P\mathrm{-algebras}).\end{gather}This equivalence states in particular that every homotopy $P$-algebra is weakly homotopy equivalent to a true $P$-algebra.

In the rest of the paper, we elaborate some important special cases. The first one is the rectification of homotopy coherent diagrams of spaces over a topological category, going back to Vogt \cite{V}, Segal \cite{Se}, and others. In fact, for an arbitrary monoidal model category $\EE$, a rectification result for $\EE$-valued homotopy coherent diagrams over an $\EE$-enriched category will be seen to be a special case of an equivalence of type (\ref{rectif}). Other examples we consider include modules over $A_\infty$-algebras, and weak maps between homotopy $P$-algebras.\vspace{2ex}

\emph{Acknowledgements:} This paper has had a slow incarnation, as the main results were already presented at the homotopy conference in Arolla in August 2004. A large part of the actual writing was done while the second author was appointed to a visiting professorship at Nice in the Spring of 2005, and he wishes to express again his gratitude to the Laboratoire Dieudonn\'e for its hospitality and support. We would also like to thank John Harper for pointing out a mistake in an earlier version of this paper.

\section{Basic definitions and examples}

Let $\EE$ be a cocomplete symmetric monoidal category.  We will assume that $\EE$ is \emph{closed}, and write $\Hom_\EE(X, Y)$ or $Y^X$ for the object of $\EE$ representing the internal $\hom$. The closedness of $\EE$ implies that the tensor product $\otimes$ of $\EE$ preserves colimits in each variable. The unit of $\EE$ will be denoted by $I$. The symmetric group on $n$ letters will be denoted by $\Sg_n$. For any finite group $\Gamma$, the category of objects in $\EE$ equipped with a right $\Gamma$-action will be denoted by $\EE^\Gamma$. It is again a cocomplete closed symmetric monoidal category.

\subsection{$C$-coloured operads}\label{def}Let $C$ be a set.  We will refer to the elements of $C$ as ``colours".  A \emph{$C$-coloured operad} $P$ is given by the following data:\vspace{1ex}

(i)  for each $n \geq 0$, and each $(n+1)$-tuple $(c_1,\dots,c_n; c)$ of colours, an object $$P(c_1,\dots,c_n;c)\text{ in }\EE;$$ 

(ii)  for each colour $c$, a unit $1_c:I \rightarrow P(c;c);$\vspace{1ex}

(iii) for each $(n+1)$-tuple $(c_1,\dots,c_n;c)$ of colours and $n$ other colour-tuples $$(d_{1,1},\dots,d_{1,{{k}_1}}),\dots,(d_{n,1},\dots,d_{n,{{k}_n}}),$$ of lengths $k_1,\dots,k_n$ respectively, a composition product\begin{gather*}P(c_1,\dots,c_n; c) \otimes P(d_{1,1},\dots,d_{1,{{k}_1}};c_1)\otimes \cdots \otimes P(d_{n,1},\dots,d_{n,k_n};c_n)\\\overset{\gamma}{\lra}P(d_{1,1},\dots,d_{n,{{k}_n}};c);\end{gather*}

(iv) for each $\sg \in \Sg_n$, a map $\sg^\star : P(c_1,\dots,c_n;c)\rightarrow P(c_{\sg(1)},\dots,c_{\sg(n)};c)$.\vspace{1ex}

The object $P(c_1,\dots,c_n;c)$ represents operations, taking $n$ inputs of colours $c_1,\dots,c_n$ respectively, and producing an output of colour $c$; this will be made precise by the definition of a coloured \emph{endomorphism-operad} below. The four data of a coloured operad are required to satisfy several axioms, which are the obvious analogues of the axioms for ordinary symmetric operads: first, the maps in (iv) define a right action by the symmetric group $\Sg_n$, in the sense that $\sg^\star\tau^\star = (\tau \sg)^\star$ for any $\sg, \tau \in \Sg_n$;  secondly, each $1_c$ is a 2-sided unit for the composition product $\gamma$;  and finally, this composition product $\gamma$ is associative and $\Sg_n$-equivariant in some natural sense.

With the obvious morphisms, the $C$-coloured operads in $\EE$ form a category, denoted $\Oper_C(\EE)$.

\subsection{$P$-algebras and coloured endomorphism-operads}\label{endo}For a $C$-coloured operad $P$, a \emph{$P$-algebra} is a family $A=(A(c))_{c \in C}$ of objects of $\EE$, together with maps $$\al_{{c}_1,\dots,c_n;c}:P(c_1,\dots,c_n;c)\otimes A(c_1)\otimes \dots \otimes A(c_n) \rightarrow A(c)$$satisfying obvious axioms for associativity, unit and equivariance.  For example, for each $\sg \in \Sg_n$, the diagram \begin{diagram}[UO,small,silent]P(c_1,\dots,c_n;c)\otimes A(c_1)\otimes \cdots \otimes A(c_n)&\rTo^{\al_{c_1,\dots,c_n;c}}&A(c)\\\dTo^{\sg^*\otimes\sg_*^{-1}}&\ruTo_{\al_{c_{\sg(1)},\dots,c_{\sg(n)};c}}&\\P(c_{\sg(1)},\dots,c_{\sg(n)};c)\otimes A(c_{\sg(1)})\otimes \cdots \otimes A(c_{\sg(n)})&&\end{diagram} commutes, where $\sg_*$ denotes the left action on tensor products induced by $\sg$ and the symmetry of $\EE$.  We will denote such an algebra by $(A,\al)$, or simply by $A$.

Equivalently, a $P$-algebra $(A,\al)$ can also be defined as a map of coloured operads $$\al:P\lra\End(A)$$ with values in the \emph{endomorphism-operad} $\End(A)$ of the family $(A(c))_{c\in C}$. This coloured operad is defined by setting$$\End(A)(c_1,\dots,c_n;c)=\Hom_\EE(A(c_1)\otimes\cdots\otimes A(c_n),A(c))$$where the composition products (resp. the $\Sg_n$-actions) are induced by substitution (resp. permutation) of the tensor factors.  

A map of $P$-algebras $f:A \rightarrow B$ is a family $(f_c:A(c)\rightarrow B(c))_{c \in \mathcal{C}}$ of maps in $\EE$, respecting the algebra structures in the sense that each diagram of the form

\begin{diagram}[UO,small]P(c_1,\dots,c_n;c)\otimes A(c_1)\otimes \cdots \otimes A(c_n) &\rTo&A(c)\\\dTo^{id \otimes f_{{c}_1}\otimes \cdots \otimes f_{{c}_n}}&&\dTo_{f_c}\\P(c_1,\dots,c_n;c)\otimes B(c_1)\otimes \cdots \otimes B(c_n)&\rTo&B(c)\end{diagram}commutes. This defines a category of $P$-algebras, denoted \emph{Alg}$_\EE(P)$. Exactly as in the uncoloured case, a map of $C$-coloured operads $\al:P\to Q$ induces adjoint functors$$\al_!:\Alg_\EE(P)\rlto\Alg_\EE(Q):\al^*.$$

\begin{rmk}\label{order}Suppose the set $C$ of colours is equipped with a linear ordering $\leq$.  If $c_1 \leq \dots\leq c_n$ are elements of $C$, write $\Sg_{c_1\dots c_n}\subseteq \Sg_n$ for the subgroup of permutations $\sg$ for which $c_{\sg (1)} \leq \dots \leq c_{\sg(n)}$ (so in particular $(c_{\sg(1)},\dots,c_{\sg(n)})$ is the same $n$-tuple as $(c_1,\dots,c_n))$.  Then a $C$-coloured operad $P$ is completely determined by the objects $P(c_1,\dots,c_n;c)$ for $c_1\leq\dots \leq c_n$, and $\Sg_{c_1\dots c_n}$ acts from the right on this object.  In other words, we can view a $C$-coloured operad as an object in $$\prod_{c_1 \leq \dots \leq c_n,c}\EE^{\Sg_{c_1\dots c_n}},$$ equipped with units $I\rightarrow P(c;c)$ and suitably equivariant and associative composition maps. Similarly, a $P$-algebra structure on a family $A=\{ A(c)\}_{c \in \mathcal{C}}$ is completely determined by ($\Sg_{c_1\dots c_n}$-equivariant) action maps$$P(c_1,\dots, c_n;c)\otimes_{\Sg_{c_1\dots c_n}}(A(c_1)\otimes \dots \otimes A(c_n))\rightarrow A(c),$$for any $c$ and any ordered sequence $c_1 \leq \dots \leq c_n$.  It will often be convenient to work with this ``smaller" representation of a coloured operad $P$ and its algebras.\end{rmk}

\begin{rmk}\label{uncoloured}If the set $C$ of colours is singleton, a $C$-coloured operad $P$ is the same as a classical (symmetric) operad, and $P$-algebra has its classical meaning. We will speak of uncoloured operads if we want to emphasize that there is just one colour. 

There is an obvious notion of \emph{non-symmetric coloured operad}, obtained from \ref{def} by deleting all references to the symmetric group actions.  If we speak of a ``$C$-coloured operad", we will always mean one in the sense of Definition \ref{def}, although we will sometimes speak of ``\emph{symmetric} $C$-coloured operads", for emphasis. Exactly as for classical operads, every non-symmetric $C$-coloured operad has an induced symmetric $C$-coloured operad which defines the same category of algebras. 
\end{rmk}

\subsection{Examples of algebras over coloured operads}

\subsubsection{Modules over operads}Let $P$ be an (uncoloured) operad, and let $C=\{ a, m\}$ be a $2$-element set.  There is a $C$-coloured operad Mod$_P$ whose algebras are pairs$$(A,M)=(A(a),A(m))$$where $A$ is a $P$-algebra and $M$ is an $A$-module:  One sets Mod$_P(c_1,\dots,c_n;c)=P(n)$ if $c=a$ and all the $c_i$ are equal to $a$ also.  And one sets Mod$_P(c_1,\dots,c_n,m)=P(n)$ when exactly one of the $c_i$ is $m$, and Mod$_P(c_1,\dots,c_n,m)=0$ in all other cases.  The structure maps of Mod$_P$ are induced by those of $P$.

\subsubsection{(Bi)modules over monoids}Write $\Ass$ for the (non-symmetric) operad whose algebras are (unitary associative) monoids in $\EE$; so $\Ass(n)=I$ for every $n\geq 0$.  There is a non-symmetric operad $LMod$ on two colours, $a$ and $m$, whose algebras are pairs $(M,E)$ where $M$ is a monoid in $\EE$ acting from the left on an object $E$ of $\EE:LMod(c_1,\dots,c_n;c)=I$ if $c=a$ and each $c_i=a$, or if $c=m$ and $c_1=\dots =c_{n-1}=a$ while $c_n=m$, and $LMod(c_1,\dots,c_n;m)=0$ in all other cases.  

There are similar operads $RMod$ on two colours and $BiMod$ on three colours, whose algebras are pairs $(M,E)$ where $M$ is a monoid acting from the right on $E$, respectively triples $(M,E,N)$ where $M$ and $N$ are monoids and $E$ is an $M$-$N$ bimodule in $\EE$.

\subsubsection{Morphisms}\label{morphisms}Let $P$ be an arbitrary (non-coloured) operad.  There is a coloured operad $P^1$ on a set $\{0,1\}$ of two colours, whose algebras are triples $(A_0,A_1,f)$ where $A_0$ and $A_1$ are $P$-algebras, and $f:A_0 \rightarrow A_1$ is a map of $P$-algebras.  Explicitly,
$$P^1(i_1,\dots,i_n;i)=\begin{cases}P(n)&\text{if }\max(i_1,\dots,i_n)\leq i;\\0&\text{otherwise}.\end{cases}$$The structure maps of $P^1$ are induced by those of $P$ (for $n=0$, we agree that $\max(i_1,\dots,i_n)=-1)$.  Given a $P^1$-algebra on two objects $A_0$ and $A_1$, the objects $P(0,\dots,0;0)$ and $P(1,\dots,1;1)$ give $A_0$, respectively $A_1$, their $P$-algebra structure; furthermore, $1:I \rightarrow P(1)$ corresponds to a map $\al:I \rightarrow P^1(0;1)$ giving a map of $P$-algebras $f:A_0\rightarrow A_1$. This coloured operad has been discussed extensively in the context of chain complexes by Markl \cite{Mar}.

\subsubsection{Categories enriched in $\EE$}\label{openriched}Let $\OO$ be a set, and consider the product $C=\OO\times \OO$.  There is a (non-symmetric) $C$-coloured operad $\Cat_\OO$ whose algebras are the $\EE$-enriched categories with $\OO$ as set of objects, and for which the maps between algebras are the functors which act by the identity on objects:  One puts$$\Cat_\OO((c_1,c_1'),\dots,(c_n,c_n');(c_0',c_{n+1}))=I$$whenever $c'_i=c_{i+1}$ for $i=0,\dots,n$, and zero in all other cases.  (In particular, for $n=0$ we have $\Cat_O(;(c,c))=I$ for each $c \in C$, providing the $\Cat_O$-algebras with the necessary identity arrows.)

\subsubsection{Diagrams in $\EE$}\label{opdiagram}Let $\CC$ be a fixed $\EE$-enriched category, with set of objects $\OO$.  There is a non-symmetric $\OO$-coloured operad $\Diag_\CC$ whose algebras are covariant $\EE$-valued diagrams on $\CC$: one puts$$\Diag_\CC(o_1,\dots,o_n;o)=\begin{cases}\Hom_\CC(o_1;o)&\text{ if }n=1;\\0&\text{ if }n>1.\end{cases}$$Composition in $\Diag_\CC$ is given by composition in $\CC$. There is of course a similar operad for contravariant diagrams.

\subsubsection{Operads}\label{opoperad}We will describe a coloured operad $S_\EE$, whose category of algebras is the category of (uncoloured) operads in $\EE$. The set of colours in this case is the set $\mathbb{N}$ of natural numbers. In fact, the operad to be defined is a coloured operad $S$ in \emph{Sets}.  Then, the strong symmetric monoidal functor$$\begin{array}{rcl}\Sets&\to&\EE\\X&\mapsto&X_\EE=\coprod_{x \in X} I\end{array}$$maps $S$ to a coloured operad $S_\EE$ whose algebras are operads in $\EE$. 

The elements of $S(n_1,\dots,n_k;n)$ are equivalence classes of triples $(T,\sg,\tau)$ where $T$ is a planar rooted tree with $n$ input edges and $k$ vertices, $\sg$ is a bijection $\{1,\dots,k\} \rightarrow V(T)$ (i.e. the set of vertices of $T$) with the property that the vertex $\sg(i)$ has valence $n_i$ (i.e. $n_i$ input edges), and $\tau$ is a bijection $\{1,\dots,n\}\rightarrow in(T)$, the set of input edges of $T$.  Two such triples $(T,\sg,\tau),(T',\sg',\tau')$ represent the same element of $S(n_1,\dots,n_k;n)$ if there is a (planar) isomorphism $\varphi:T\rightarrow T'$ with $\varphi \circ \tau=\tau'$ and $\varphi \circ \sg=\sg'$.

Any $\al \in \Sg_k$ induces a map $\al^\star:S(n_1,\dots,n_k;n)\rightarrow S(n_{\al(1)},\dots,n_{\al(k)};n)$ sending (the equivalence class of) $(T, \sg, \tau)$ to $(T, \sg \al, \tau)$.  The identity element $1_n \in S(n;n)$ is represented by the tree $t_n$ (the corolla with $n$ leaves) whose inputs are numbered $1,\dots,n$ from left to right with respect to the planar structure.

The composition product is defined as follows: given $(T,\sg,\tau)$ as above, and $k$ other such $(T_1,\sg_1,\tau_1),\dots,(T_k,\sg_k,\tau_k)$, with $n_1,\dots,n_k$ inputs and $p_1,\dots,p_k$ vertices respectively, one obtains a new planar rooted tree $T'$ by replacing the vertex $\sg(i)$ in $T$ by the tree $T_i$, identifying the $n_i$ input edges of $\sg(i)$ in $T$ with the $n_i$ input edges of $T_i$ via the bijection $\tau_i$ (the $l$-th input edge of $\sg(i)$ in the planar order is matched with the input edge $\tau_i(l)$ of $T_i$).

The vertices of the new tree $T'$ are numbered in the following order: first the vertices of $T_{\sg(1)}$ in the order given by $\sg_1$, then the vertices in $T_{\sg(2)}$ in the order given by $\sg_2$, etc.  In other words, the map
$\{1,\dots, p_1+\dots +p_k\}\rightarrow V(T')$ is given by $(\sg_1\times\cdots\times\sg_k)\circ\sg(p_1,\dots,p_k)$ where $\sg(p_1,\dots,p_k)$ permutes the blocks of size $p_i$. The new tree $T'$ still has $n$ input edges, which are ordered as given by $\tau$ and the identifications given by the $\tau_i$.  Notice that $S(n_1;n)=\Sg_n$ if $n_1=n$, and $S(n_i,n)=\phi$ otherwise.  More precisely, $S(n,n)$ consists of pairs $(t_n,\tau)$ where $t_n$ is the tree above and $\tau$ is a numbering of its inputs.  The composition product of $S$ in particular gives a map $S(n,n)\times S(n,n)\rightarrow S(n,n)$ which sends $((t_n,\tau),(t_n,\rho))$ to $(t_n,\rho\tau)$, so that $S(n;n)$ is identified with the \emph{opposite group} of $\Sg_n$.  

The $S$-algebras are exactly the operads in sets. Indeed, given such an operad $P$, a triple $(T,\sg,\tau)\in S$ acts on $(p_1,\dots,p_k)\in P(n_1)\times\cdots\times P(n_k)$ by labelling the vertex $\sg(i)\in T$ by $p_i$, and then using the operad structure of $P$ to compose $(p_1,\dots,p_k)$ along the tree $T$ to get an element $p\in P(n)$, and then applying the right action by $\tau$ to this element; i.e.$$(T,\sg,\tau)(p_1,\dots,p_k)=p\cdot \tau.$$In particular, the action $S(n,n)\times P(n) \rightarrow P(n)$ encodes the right $\Sg_n$-action.

There is a similar coloured operad $S^+$, whose algebras are operads $P=(P(n))_{n\geq 1}$ without $0$-term. It is the coloured suboperad of $S$ given by considering only vertices of valence $\geq 1$. This operad $S^+$ is in fact induced by a non-symmetric $\mathbb{N}$-coloured operad, cf. Remark \ref{uncoloured}. Indeed, it is sufficient to show that for any planar tree $T$ and any ordering $\tau$ of its input edges, there is an ordering $\sg_{T,\tau}$ of its vertices, which is compatible with the composition product just described, and is invariant in the sense that for a planar isomorphism $\varphi:T\rightarrow T',$ the following relation holds:$$\sg_{T',\varphi \circ \tau}=\varphi \circ \sg_{{T,}\tau}.$$The easiest way to define $\sg_{T,\tau}$ is to view $\sg_{T,\tau}$ (resp. $\tau$) as linear orderings of the vertices (resp. input edges) of $T$. So, given a linear ordering of the input edges of $T$, we have to define a linear ordering of its vertices; we use induction on $T$.  

Suppose that $T$ is obtained by grafting $p$ trees on the corolla $t_p$, for short: $T=t_p(T_1,\dots,T_p)$, with root $r$; write $\leq_i$ for the linear order on $T_i$ corresponding to the linear ordering of its input edges induced by that of $T$.  Also write $T_i<T_j$ if the first input edge (in the ordering) of $T_i$ comes before the first one of $T_j$. This defines a linear order on the $p$-element set $\{T_1,\dots,T_p\}$.  For vertices $v, w \in T$, now define$$v \leq w \iff\begin{cases}&v = r;\\\text{or }&v,w \in T_i\text{ and }v \leq_i w;\\\text{or }&v\in T_i\text{ and }w \in T_j\text{ and }T_i < T_j.\end{cases}$$

\subsubsection{Coloured operads}\label{opcoloperad}Let $C$ be a set of colours. There is a coloured operad $S_C$ whose algebras are $C$-coloured operads. The set of colours of $S_C$ is the set of sequences of the form $(c_1,\dots,c_n;c)$ for all $n \geq 0$ (or, more economically, those sequences for which $c_1 \leq \dots \leq c_n$ after having chosen an order on $C)$.  We leave a detailed description of $\SS_C$ to the reader.

\subsection{Change of colour}\label{changecol}If $P$ is a $C$-coloured operad and $\al:D\to C$ is a map between sets ``of colours'', then $P$ pulls back to a $D$-coloured operad $\al^*(P)$ in the obvious way,$$\al^*(P)(d_1,\dots,d_n;d)=P(\al(d_1),\dots,\al(d_n);\al(d)).$$This defines a functor$$\al^*:\Oper_C(\EE)\to\Oper_D(\EE).$$In this way, $C$-coloured operads for varying sets of colours $C$ together form a fibered category over the category of sets. Objects of this fibered category are pairs $(C,P)$ where $P$ is a $C$-coloured operad, and arrows $(D,Q)\to(C,P)$ are pairs $(\al,\phi)$ where $\al:D\to C$ and $\phi:Q\to\al^*(P)$ is a map of $D$-coloured operads. Such an arrow induces an adjoint pair $$(\al,\phi)_!:\Alg_\EE(Q)\rlto\Alg_\EE(P):(\al,\phi)^*.$$For instance, if $S_C$ and $S_D$ are the operads whose algebras are $C$-coloured operads and $D$-coloured operads respectively, the map $\al:D\to C$ induces in this way a map of coloured operads $S_D\to S_C$, and hence an adjoint pair$$\al_!:\Oper_D(\EE)\rlto\Oper_C(\EE):\al^*.$$

\subsubsection{Example}Let $P$ be an uncoloured operad. A graded $P$-algebra is a sequence $(A_n)_{n\geq 0}$ of objects of $\EE$ (indexed by $n\in\NN$), such that the coproduct $A=\coprod A_n$ has a $P$-algebra structure, which respects the grading in the sense that the structure map $P(k)\otimes A^{\otimes k}\to A$ maps the summand $P(k)\otimes A_{k_1}\otimes\cdots\otimes A_{n_k}$ to the summand $A_n$ for $n=n_1+\cdots+n_k$. Associated to $P$, there is an $\NN$-coloured operad $\Gr(P)$ whose algebras are the graded $P$-algebras. It is given by$$\Gr(P)(n_1,\dots,n_k;n)=\begin{cases}P(k)&\text{if }n_1+\cdots+n_k=n;\\0&\text{otherwise}.\end{cases}$$There is an evident map of coloured operads$$(\al,\phi):(\NN,\Gr(P))\lra(*,P)$$given by the unique map $\al:\NN\to *$ and the inclusion $\phi:\Gr(P)(n_1,\dots,n_k;n)\ito\al^*(P)(n_1,\dots,n_k;n)=P(k)$. The left adjoint $(\al,\phi)_!$ sends the graded $P$-algebra $(A_n)_{n\geq 0}$ to the $P$-algebra $\coprod A_n$.

\subsubsection{Example}Suppose $\EE$ is additive. Let $\PLie$ be the operad for pre-Lie algebras in $\EE$, and let $S_0$ be the $\NN$-coloured operad for non-symmetric operads in $\EE$. ($S_0$ is defined as $S$ in \ref{opcoloperad}, but without the $\tau$'s). There is a map of $\NN$-coloured operads$$(\sg;c):(\NN,\Gr(\PLie))\to(\NN,S_0)$$defined as follows: $\sg(n)=n+1$, and $c$ sends the pre-Lie operation $\oo_{n,m}\in\Gr(\PLie)(n,m;n+m)$ to the sum of the $\oo_i$-operations in $S_0(n+1,m+1;n+m+1)$, where the $\oo_i$-operation is the tree with two vertices of valence $n+1$ and $m+1$ respectively and one internal edge at the $i$-th entry of the lower vertex. The familiar construction of a pre-Lie algebra out of a non-symmetric operad is now given by the composition $(\al,\phi)_!\circ(\sg,c)^*$, where\begin{diagram}[UO,small](*,\PLie)&\lTo^{(\al,\phi)}&(\NN,\Gr(\PLie))&\rTo^{(\sg,c)}&(\NN,S_0)\end{diagram}as above. (There is also a symmetric version of this construction, cf. Kapranov-Manin \cite{KM}.)

The functorial Boardman-Vogt resolution to be discussed in Section $4$ will evidently have the property that $\Gr(\WW(P))=\WW(\Gr(P))$, and the maps above will induce maps between the corresponding cofibrant resolutions\begin{diagram}[UO,small](*,\WW(\PLie))&\lTo&(\NN,\WW(\Gr(\PLie)))&\rTo&(\NN,\WW(S_0)).\end{diagram}This shows that the same construction yields for any ``operad up to homotopy'' a pre-Lie algebra up to homotopy (and hence a $\WW(\mathrm{Lie})$-algebra, i.e. an $L_\infty$-algebra).

\section{Model structure on $P$-algebras}In this section, we assume that our cocomplete symmetric monoidal closed category $\EE$ is equipped with a compatible Quillen model structure, making it into a so-called \emph{monoidal model category}.  We will always assume that the unit $I$ of $\EE$ is cofibrant and that the model structure is cofibrantly generated.  Recall that under the last assumption, for any finite group $G$, there is an induced monoidal model structure on the category $\EE^G$ of objects with right $G$-action, for which the forgetful functor $\EE^G \rightarrow \EE$ preserves and reflects fibrations and weak equivalences.  We refer to \cite{Hir}, \cite{Hov} for basic facts concerning monoidal model categories and associated equivariant categories like $\EE^G$, see also \cite[2.5]{BM2}.

For a set of colours $C$ and a $C$-coloured operad $P$, our purpose is to show that, under suitable conditions, the category $Alg_\EE(P)$ of $P$-algebras admits a closed model structure for which the forgetful functor$$U_P:\Alg_\EE(P)\rightarrow \EE^{C}$$preserves and detects fibrations and weak equivalences.  In other words, a map of $P$-algebras $A \rightarrow B$ is a fibration (resp. weak equivalence) in $\Alg_\EE(P)$ if and only if for each colour $c \in C$, the map $A(c)\rightarrow B(c)$ is a fibration (respectively weak equivalence) in $\EE$.  Let us call a coloured operad \emph{admissible} if this defines a closed model structure on $Alg_\EE(P)$.  This terminology extends the one in \cite{BM1} for uncoloured operads.  Our first result is an immediate generalisation of a result from \cite{BM1} for uncoloured operads.  Recall that a (cocommutative) coalgebra interval $H$ is a (cocommutative) counital comonoid object $H$ together with a factorisation of the codiagonal $\nabla: I\sqcup I \rightarrow I$ (which is a map of comonoids) as$$I\sqcup I \overset{(i_0,i_1)}{\ito}H \overset{\epsilon} {\lra}I$$where both maps are maps of comonoids, while $(i_0,i_1)$ is a cofibration and $\epsilon$ is a weak equivalence in $\EE$. In particular, the monoidal model categories of compactly generated spaces, of simplicial sets and of symmetric spectra admit such cocommutative coalgebra intervals, while the category of chain complexes admits a coalgebra interval.

\begin{thm}\label{main0}Let $\EE$ be a (cofibrantly generated) monoidal model category with cofibrant unit and a symmetric monoidal fibrant replacement functor. If $\EE$ has a coalgebra interval, then every non-symmetric coloured operad is admissible. If the interval is moreover cocommutative, the same is true for every symmetric coloured operad.\end{thm}
\begin{proof}The proof is identical to the proof of \cite[Theorem 3.2]{BM1} and the remark following it.  The crucial observation is that, if $A$ is a fibrant $P$-algebra, the family $A^H=\{A(c)^H\}_{c \in C}$ with the induced maps$$A\stackrel{\epsilon^\star}{\rightarrow} A^H \stackrel{i_0^\star,i_1^\star}{\rightarrow} A \times A$$provides a path-object for $A$ in the category of $P$-algebras.\end{proof}

\section{Free operads and the Boardman-Vogt resolution}

Let $\EE$ be a monoidal model category, always assumed to have a cofibrant unit, as before. Let $C$ be a set of colours. Again it will be convenient to choose a linear order on $C$.

The category of $C$-coloured operads is itself a category of algebras over a coloured (non-symmetric) operad, cf. Example \ref{opoperad}. Thus Theorem \ref{main0} provides conditions under which this category carries a model structure, in which a map $Q\to P$ between $C$-coloured operads is a fibration (respectively, a weak equivalence) if and only if, for each sequence of colours $c_1,\dots,c_n,c,$ the map $Q(c_1,\dots,c_n;c)\to P(c_1,\dots,c_n;c)$ in $\EE$ is a fibration (respectively, a weak equivalence). Even, when this does not give a model structure on the category of $C$-coloured operads, one still can call a map a \emph{cofibration} when it has the left lifting property with respect to these trivial fibrations, and an object $R$ \emph{cofibrant} if the unique map $I_C\to R$ is a cofibration. Here $I_C$ is the initial $C$-coloured operad, $I_C(c;c)=I$, and all other components of $I_C$ are zero. Throughout this section and the next, we will use the terminology cofibrant/cofibration in this sense.

\subsection{Coloured collections}--\\
 
A $C$-collection $K$ is given by objects $K(c_1,\dots,c_n;c)$ and right actions $$\sg^*:K(c_1,\dots,c_n;c)\to K(c_{\sg(1)},\dots,c_{\sg(n)};c)$$for each $\sg\in\Sg_n$, exactly as for $C$-coloured operads in Definition \ref{def}. In particular, if $c_1\leq\cdots\leq c_n$ then $\Sg_{c_1\dots c_n}$ acts from the right on $K(c_1,\dots,c_n;c)$. We call $K$ a \emph{$\Sg$-cofibrant $C$-collection} if each such $K(c_1,\dots,c_n;c)$ is cofibrant in $\EE^{\Sg_{c_1\dots c_n}}$. The $C$-collection $K$ is called \emph{pointed} if it is equipped with units $1_c:I\to K(c;c)$, one for each colour $c\in C$. A pointed collection $K$ is \emph{well-pointed} if these units are cofibrations in $\EE$. With the obvious notion of maps, this defines categories $\Coll_C(\EE)$ of $C$-coloured collections, and $\Coll_C^\star(\EE)$ of pointed $C$-coloured collections. Notice that both categories carry a model structure for which a map $K\to L$ is a weak equivalence (resp. a fibration) if and only if for each sequence of colours $c_1,\dots,c_n,c$, the map $K(c_1,\dots,c_n;c)\to L(c_1,\dots,c_n;c)$ is a weak equivalence (resp. a fibration) in $\EE$. The $\Sg$-cofibrant collections are the cofibrant objects in $\Coll_C(\EE)$, while the well-pointed $\Sg$-cofibrant collections are the cofibrant objects in $\Coll_C^\star(\EE)$. 

\begin{thm}\label{forget}The forgetful functor from $C$-coloured operads to pointed $C$-coloured collections has a left adjoint $F_C^\star:\Coll_C^\star(\EE)\to\Oper_C(\EE).$\end{thm}

We will give an explicit construction in a moment. But notice that this is really a special case of the adjunction $\phi_!:\Alg_\EE(P)\lrto\Alg_\EE(Q):\phi^*$, considered in Section \ref{endo}, because both categories $\Coll_C^\star(\EE)$ and $\Oper_C(\EE)$ are categories of algebras for suitable coloured operads, cf. Example \ref{opcoloperad}.

\begin{cor}The free $C$-coloured operad on a well-pointed $\Sg$-cofibrant $C$-collection is a cofibrant $C$-coloured operad.\end{cor}

Recall that this has a sense even when there is no model structure on the category of $C$-coloured operads. The corollary follows from the theorem by adjunction.

It is useful to give an explicit description of the free $C$-coloured operad $F^\star_C(K)$ on a pointed $C$-collection, exactly as in the uncoloured case. For this, recall the groupoid $\TT$ of planar rooted trees and non-planar isomorphisms as in \cite{BM1,BM2}. There is a similar groupoid $\TT_C$ of planar rooted trees all of whose edges are coloured by elements of $C$. The planar structure of the tree induces a linear order on the input edges of any given vertex; it will be convenient to require that these edges are coloured in an order-preserving way. The arrows in $\TT_C$ are (non-planar) isomorphisms which preserve the colours of the edges. For each sequence $c_1\leq\cdots\leq c_n$ of colours and each colour $c$, there is a special object of $\TT_C$, the corolla tree $t_n(c_1,\dots,c_n;c)$, with only one vertex, and $n$ input edges coloured $c_1,\dots,c_n$ respectively, and one output edge coloured $c$. The automorphism group of this tree is exactly $\Sg_{c_1\dots c_n}$.

For each $n\geq 0$, let $\TT_C(n)$ be the full subcategory of $\TT_C$ given by the trees with $n$ input edges. For a coloured tree $T$ in $\TT_C(n)$, write $\inp(T)$ for the set of these input edges, and $\lda_n(T)$ for the set of all bijections $\{1,\dots,n\}\to\inp(T).$ Then $\Sg_n$ acts from the right on $\lda_n(T)$, so $\lda_n$ defines a functor $\lda_n:\TT_C(n)\to\Sets^{\Sg_n}$. The planar structure of a tree $T$ induces an order on the input edges of $T$ which allows us to identify $\lda_n(T)$ with $\Sg_n$ itself. The functoriality of $\lda$ then gives a left action of $\Aut(T)$ on $\Sg_n$ as a right $\Sg_n$-set (not as a group).

Intuitively speaking, elements of the free operad on the pointed collection $K$ are represented by trees with\begin{itemize}\item[-]inputs labelled by $1,\dots,n$;\item[-]edges ``coloured'' by colours of $C$;\item[-]vertices labelled by elements of $K$,\end{itemize}in a compatible way so that a vertex is labelled by an element $k\in K(c_1,\dots,c_n;c)$ only if the incoming edges are labelled (from left to right) by the colours $c_1,\dots,c_n$, and the outgoing edge is labelled by $c$. Furthermore, some identifications are made, taking account for tree-automorphisms and for identities.

For a formal definition, we define for each $C$-coloured pointed collection $K$ a functor$$\Ku:\TT_C^\op\to\EE$$by induction on trees, exactly as in \cite[5.8]{BM1}:\begin{align*}\Ku(|_c)&=\KK(;c),\\\Ku(t_n(c_1,\dots,c_n;c))&=\KK(c_1,\dots,c_n;c),\\\Ku(T)&=\KK(c_1,\dots,c_n;c)\otimes\Ku(T_1)\otimes\cdots\otimes\Ku(T_n),\end{align*}where $T=t_n(c_1,\dots,c_n;c)[T_1,\dots,T_n]$ is obtained by grafting the trees $T_i$ with output colour $c_i$ on the top of the $i$-th input edge of the coloured tree $t_n(c_1,\dots,c_n;c)$. The free operad $F^\star_C(\KK)$ is now obtained by first constructing the free operad on an unpointed collection$$F_C(\KK)(n)=\coprod_{[T],T\in\TT_C(n)}\Ku(T)\otimes_{\Aut(T)}I[\Sg_n],$$and then factoring out the identities by constructing the pushout square of coloured operads\begin{diagram}[UO,small,silent]F(I_C)&\rTo&F_C(\KK)\\\dTo&&\dTo\\I_C&\rTo&F_C^\star(\KK).\end{diagram}

\subsection{Segments and intervals}Recall from \cite[4.1-2]{BM2} that a \emph{segment} in $\EE$ is an object $H$, equipped with maps\begin{diagram}[UO,small]I&\pile{\rTo^0\\\rTo_1}&H&\rTo^\eps&I\end{diagram}and an associative operation $\vee:H\otimes H\to H$, having $0$ as a unit element, $1$ as an absorbing element, and $\eps$ as a counit. (Set-theoretically, these axioms on $H$ can be paraphrased by the equations $0\vee x=x=x\vee 0,\,1\vee x=1=x\vee 1,\,\eps(x\vee y)=\eps(x)\eps(y),\,\eps\circ 0=id_I=\eps\circ 1$.) An \emph{interval} is such a segment for which $I\sqcup I\overset{(0,1)}{\lra}H$ is a cofibration and $H\overset{\eps}{\lra}I$ is a weak equivalence. Notice that $I\sqcup I$ and $I$ are segments, but (usually) not intervals. The standard example is the real unit interval $[0,1]$ with its maximum operation, in the model category $\Top$ of compactly generated topological spaces. But there are many other examples, cf. loc. cit.

An interval $H$ allows us to construct for each $C$-coloured operad $P$ a new operad $\WW(H,P)$, the Boardman-Vogt resolution of $P$. This operad is constructed like the free operad $F^\star_C(P)$, but with the \emph{additional} assignment of a ``length'' in $H$ to each internal edge of the trees. The axioms for segments enter into the identifications to be made in the construction of $\WW(H,P)$, which, besides the identifications coming from automorphisms of trees, are of the following two types:

\begin{itemize}\item[(i)]edges of length $0\in H$ are contracted, using the operad structure of $P$;\item[(ii)]edges around a vertex labelled by a unit $1_C\in P(c;c)$ are contracted into a single edge, deleting the vertex and assigning the sup ($\vee$) of the corresponding lengths as new length.\end{itemize}The properties of this construction are summarised in the following theorem. Here, we call a $C$-coloured operad \emph{$\Sg$-cofibrant} (resp. \emph{well-pointed}) if the forgetful functor of Theorem \ref{forget} maps it to a $\Sg$-cofibrant (resp. well-pointed) collection. The proof proceeds along exactly the same lines as the one for the uncoloured case, described in detail in \cite{BM2}.

\begin{thm}\label{W}Let $\EE$ be a (cofibrantly generated) monoidal model category with cofibrant unit $I$ and interval $H$.

(i) For each well-pointed $\Sg$-cofibrant $C$-coloured operad $P$, the counit of the free-forgetful adjunction factors into a cofibration followed by a weak equivalence$$F^\star_C(P)\ito\WW(H,P)\eqv P$$of $C$-coloured operads. This factorisation is natural in $P$ and $H$.

(ii) If $P\to Q$ is a $\Sg$-cofibration between well-pointed $\Sg$-cofibrant $C$-coloured operads, then the induced map $\WW(H,P)\to\WW(H,Q)$ is a cofibration of cofibrant $C$-coloured operads.\end{thm}

As in \cite{BM2} there is also a \emph{relative version} of the coloured Boardman-Vogt resolution. It provides for each $\Sg$-cofibration $u:P\to Q$ of $C$-coloured operads a factorisation of $u$ into a cofibration followed by a weak equivalence$$P\ito\WW(H,P\overset{u}{\lra}Q)\eqv Q.$$In fact, it has a somewhat stronger property: let $P[u]$ be the free extension of $P$ by $u$, constructed as a  pushout\begin{diagram}[UO,small,silent]F^\star_C(P)&\rTo&F^\star_C(Q)\\\dTo&&\dTo\\P&\rTo&P[u]\end{diagram}of $C$-coloured operads. Then $P\to P[u]$ is a cofibration of $C$-coloured operads, and the cofibration $P\ito\WW(H,P\to Q)$ factors as two cofibrations$$P\ito P[u]\ito\WW(H,P\overset{u}{\lra}Q).$$The operad $\WW(H,P)$ describes ``$P$-algebras up to homotopy''. For example, in the case of topological spaces with the standard unit-interval, or the case of chain complexes with normalized chains on the standard $1$-simplex as interval, one can take the operad $\Ass$ for associative monoids and finds a subdivision of Stasheff's $A_\infty$-operad as $\WW(H,\Ass)$. Stasheff's $A_\infty$-operad itself can be obtained as the relative Boardman-Vogt resolution $\WW(H, I_*\to \Ass)$ where $I_*$ is the operad for pointed objects. The relative coloured case enables us also to consider, for example, an operad whose algebras are pairs $(R,M)$ where $R$ is a (strict) ring, but $M$ is an $R$-module only up to homotopy.

\subsection{Change of colour}\label{changecol2}Recall from Section \ref{changecol} for a map $\al:D\to C$ between sets of colours the adjunction $\al_!:\Oper_D(\EE)\lrto\Oper_C(\EE):\al^*$ between $D$-coloured operads and $C$-coloured operads. For a $D$-coloured operad $Q$ and a $C$-coloured operad $P$, there are natural maps$$\al_!\WW(H,Q)\to\WW(H,\al_!Q)\text{ and }\WW(H,\al^*P)\to\al^*\WW(H,P),$$but in general these maps are not isomorphisms. If $\al$ is one-to-one, however, there is an explicit description of $\al_!(P)$, as$$\al_!(P)(d_1,\dots,d_n;d)=\begin{cases}P(c_1,\dots,c_n;c)&\text{if }d_i=\al(c_i),d=\al(c),\\I&\text{if }n=1,d=d_1\not\in Im(\al),\\0&\text{otherwise},\end{cases}$$and in this case, the map $\al_!\WW(H,Q)\to\WW(H,\al_!Q)$ is an isomorphism.

\section{The comparison theorem}

As before, we work with a monoidal model category $\EE$, and fix a set of colours $C$. We call a map $\varphi:P \rightarrow Q$ between $C$-coloured operads a \emph{weak equivalence} if each of its components $\varphi_{c_1,\dots c_n;c}:P(c_1,\dots,c_n;c)\to Q(c_1,\dots,c_n;c)$ is a weak equivalence in $\EE$, independently of whether this is part of a model structure or not.
 
As already pointed out in Section \ref{endo}, a map $\varphi:P \rightarrow Q$ of $C$-coloured operads induces an evident functor$$\varphi^\star: Alg_\EE(Q)\rightarrow Alg_\EE(P),$$commuting with the forgetful functors, i.e. $U_P \circ \varphi^\star=U_Q.$  For general reasons, this functor $\varphi^\star$ has a left adjoint, which we denote by$$\varphi_!:Alg_\EE(P) \rightarrow  Alg_\EE(Q).$$This functor sends free $P$-algebras to free $Q$-algebras, i.e.$$\varphi_!(F_P(X))=F_Q(X),$$and this can in fact be used to give an explicit description of $\varphi_!$: Given a $P$-algebra $A$, we can write $A$ as a coequaliser of free $P$-algebras, \begin{diagram}[UO,small]F_P U_P F_P U_P(A)&\pile{\rTo^\beta\\\rTo_\gamma}&F_P U_P (A)&\rTo^\al A\end{diagram}where $\al$ is the $P$-algebra structure on $A$, $\beta$ is $F_PU_P(\al)$, and $\gamma$ is the $P$-algebra structure on $F_PU_P(A)$.  Then $\varphi_!(A)$ is the coequaliser \begin{diagram}F_QU_PF_PU_P(A)&\pile{\rTo^{\beta'}\\\rTo_{\gamma'}}&F_QU_P(A)&\rTo\varphi_!(A)\end{diagram}where $\beta'=F_QU_P(\al)$, and $\gamma'$ is the unique $Q$-algebra map extending the map $F_PU_P(A)\to F_QU_P(A)$ in $\EE^C$ given by $\varphi$.

\begin{thm}\label{comparison}Let $\EE$ be a left proper (cofibrantly generated) monoidal model category with cofibrant unit and $\varphi:  P \rightarrow Q$ a weak equivalence between admissible $\Sg$-cofibrant well-pointed $C$-coloured operads in $\EE$.  Then the adjunction $\varphi_!:\Alg_\EE(P)\lrto\Alg_\EE(Q):\varphi^\star$ is a Quillen equivalence.\end{thm}

\begin{proof}The proof follows exactly the same pattern as that of Theorem 4.4 in \cite{BM1}, see also the Appendix of \cite{BM3}, and we only give an outline. Evidently, $\phi^*$ preserves weak equivalences and fibrations, so by adjunction $\phi_!$ preserves cofibrations and trivial cofibrations, and $(\phi_!,\phi^*)$ is a Quillen pair. Since $\phi^*$ also reflects weak equivalences, it suffices \cite[1.3.16]{Hov} to prove that for each \emph{cofibrant} $P$-algebra $A$, the unit of the adjunction $\eta_A:A\to\phi^*\phi_!(A)$ is a weak equivalence. To this end, let us first recall some terminology. For a $P$-algebra $A$ and a cofibration $u:U_P(A)\ito Z$ in $\EE^C$, the \emph{free extension} of $A$ by $u$, denoted $A\ito A[u]$, is by definition the pushout in the category of $P$-algebras of $F_P(u)$ along the counit $\eps$:\begin{diagram}[small,UO,silent]F_PU_P(A)&\rTo^\eps&A\\\dTo&&\dTo\\F_P(Z)&\rTo&A[u]\end{diagram}A cellular extension of a $P$-algebra $A$ is a map $A\ito B$ which can be written as a sequential colimit of a sequence$$A=A_0\ito A_1\ito A_2\ito\cdots\ito A_\xi\ito A_{\xi+1}\ito$$indexed by some ordinal, where $A_\xi\ito A_{\xi+1}$ is such a free extension by a cofibration, and $A_\lambda=\lim_{\xi<\lambda}A_\xi$ at limit ordinals. A \emph{cellular $P$-algebra} is a cellular extension of the initial $P$-algebra $P(0)$. Since every cofibrant $P$-algebra is a retract of a cellular $P$-algebra $A$, it now suffices to prove that the unit $A\to\phi^*\phi_!(A)$ is a weak equivalence for each cellular $P$-algebra $A$. We do this by proving a stronger statement, by induction on $A$, as follows. For each cellular $P$-algebra $A$ we construct a $\Sg$-cofibrant well-pointed operad $P_A$ whose algebras are $P$-algebras under $A$, and a similar operad $Q_{\phi_!(A)}$. We then prove by induction on the cellular $P$-algebra $A$ that the map $P\to Q$ induces a weak equivalence $P_A\to Q_{\phi_!(A)}$. Since in degree $0$, this is the map $P_A(0)=A\to\phi^*\phi_!(A)=\phi^*(Q_{\phi_!A})$ (when viewed as a map of $P$-algebras), this yields the desired conclusion.

As to the induction, for the initial step $A=P(0)$ one has $P_A=P$ and $Q_{\phi_!(A)}=Q_{Q(0)}=Q$, and hence $P_A\to Q_{\phi_!(A)}$ is a weak equivalence (between $\Sg$-cofibrant well-pointed operads) by assumption. At a limit stage $A_\lambda=\lim_{\xi<\lambda}A_\xi$, we have $P_A=\lim P_{A_\xi}$ and $Q_{\phi_!A}=\lim Q_{\phi_!(A_\xi)}$, and the conclusion follows from the fact that a colimit of a ladder of weak equivalences between cofibrant objects in $\EE^C$ is again a weak equivalence, cf. \cite[2.3]{BM2}. The difficult part is the case where we have constructed the weak equivalence $P_A\to Q_{\phi_!(A)}$, and want to deduce a similar weak equivalence $P_{A[u]}\to Q_{\phi_!(A[u])}$ for a free extension $A[u]$ of $A$. By replacing $P$ by $P_A$ and $Q$ by $Q_{\phi_!(A)}$, we may assume that $A=P(0)$ and that the free extension is by a cofibration $u:U_P(P(0))\ito Z$. Then $\phi_!(A)=Q(0)$, so that $Q_{\phi_!(A[u])}=Q_{Q(0)[v]}$ where $v$ is the pushout of $u$ along $U_P(P(0))\eqv U_Q(Q(0))$ in $\EE^C$:

\begin{diagram}[UO,small,silent]U_P(P(0))&\rTo^{\sim}& U_Q(Q(0))\\\dTo^u&&\dTo_v\\Z&\rTo^{\sim}&W\end{diagram}

The map $Z\to W$ is again a weak equivalence because $\EE$ (and hence $\EE^C$) is left proper, by assumption. Now in this case, there is an explicit construction of $P_{A[u]}$, which is completely analogous to the one in \cite[5.11]{BM1} for uncoloured operads, except that one has to work with trees whose edges are labelled by elements of the fixed set $C$ of colours, as in the description of the free operad in the previous section. This construction, applied to $P$ and $u$ as well as to $Q$ and $v$, shows that $P_{P(0)[u]}\to Q_{Q(0)[v]}$ is a weak equivalence of $\Sg$-cofibrant well-pointed operads, exactly as in \cite[5.12]{BM1}.\end{proof}

\section{Rectification of diagrams}

As a first application, we will prove that any homotopy coherent diagram is equivalent to a strict diagram. For the categories of topological spaces and of simplicial sets, results of this type have received considerable attention, and go back to Vogt \cite{V}, see also Bousfield-Kan \cite{BK}, Segal \cite{Se}, Cordier-Porter \cite{CP}, Dwyer-Kan-Smith \cite{DKS}, and others. The purpose of this section is to show that rectification results of this kind, expressed more explicitly in terms of a Quillen equivalence, can be obtained in the more general context of monoidal model categories, as an immediate consequence of our Comparison Theorem \ref{comparison}.

By way of example, let us consider diagrams with values in the category $\Top$ of compactly generated topological spaces. Let $\CC$ be a small category. A \emph{$\CC$-diagram of spaces} is a functor $X:\CC\to\Top$. A \emph{homotopy coherent $\CC$-diagram} is similarly given by a space $X(c)$ for each object $c$ of $\CC$, and a map $X(\al):X(a)\to X(b)$ for each arrow $\al:a\to b$ in $\CC$. These are required to be functorial up to a given homotopy $H_{\be,\al}$, connecting $X(\be)\circ X(\al)$ and $X(\be\al)$, for any $2$-simplex $a\overset{\al}{\lra}b\overset{\be}{\lra}c$ of arrows in $\CC$. Furthermore, for any $3$-simplex of arrows $a\overset{\al}{\lra}b\overset{\be}{\lra}c\overset{\gamma}{\lra}d$, there has to be a connecting higher homotopy $H_{\gamma,\be,\al}$ between $(H_{\gamma,\be}\star X(\al))\circ H_{\gamma\be,\al}$ and $(X(\gamma)\star H_{\be,\al})\circ H_{\gamma,\be\al}$, where $\star$ denotes horizontal composition of homotopies. These $H_{\gamma,\be,\al}$ are then required to satisfy a coherence condition for any $4$-simplex of arrows in $\CC$, and so on.

Diagrams of spaces on $\CC$ are \emph{precisely} the algebras for the coloured operad $\Diag_\CC$ of Example \ref{opdiagram}, where the ambient category $\EE$ is $\Top$ and $\CC$ is viewed as a discrete topological category. Let $H=[0,1]$ be the usual interval, with binary operation $x\vee y=\max(x,y)$. Then, as pointed out by Vogt \cite{V} and Cordier-Porter \cite{CP}, homotopy coherent $\CC$-diagrams are \emph{precisely} the algebras for the coloured Boardman-Vogt resolution $\WW(H,\Diag_\CC)$. This coloured operad has unary operations only, just like $\Diag_\CC$. For objects $a$ and $b$ in $\CC$, an operation in $\WW(H,\Diag_\CC)(a;b)$ is a string of composable arrows in $\CC$,\begin{diagram}[UO,small]a&\rTo^{\al_0}&c_1&\rTo^{\al_1}&c_2&\rTo&\cdots&\rTo&c_{n-1}&\rTo^{\al_{n-1}}&c_n&\rTo^{\al_n}&b\end{diagram}together with ``waiting times''$$t_i\in[0,1],\,i=1,\dots,n,$$to be thought of as assigned to the $c_i$. The action of such a string on a point $x\in X(a)$ can be thought of as acting by $\al_0$, then waiting for $t_1$ seconds before acting by $\al_1$, then waiting for $t_2$ seconds before acting by $\al_2$, etc. Such an operation is identified with the operation\begin{diagram}[UO,small]a&\rTo^{\al_0}&c_1&\cdots&\rTo&c_{i-1}&\rTo^{\al_i\al_{i-1}}&c_{i+1}&\rTo&\cdots&c_n&\rTo^{\al_n}&b\end{diagram}and waiting times $(t_1,\dots,\hat{t}_i,\dots,t_n)$ whenever $t_i=0$, and with\begin{diagram}[UO,small]a&\rTo^{\al_0}&c_1&\cdots&\rTo&c_i=c_{i+1}&\rTo&\cdots&c_n&\rTo^{\al_n}&b\end{diagram} and waiting times $(t_1,\dots,t_i\vee t_{i+1},\dots,t_n)$ whenever $\al_i=id_{c_i}$. Composition of operations is defined as concatenation of strings, assigning waiting time $1$ to the intermediate arrow.

There is an evident map of coloured operads $\WW(H,\Diag_\CC)\overset{\eps}{\lra}\Diag_\CC$, which forgets the waiting times and composes the string. Pulling back along this map $\eps$ is the functor which views each strict diagram as a homotopy coherent diagram with trivial homotopies.

Now consider $\Top$ as a Quillen model category in the standard way with weak homotopy equivalences and Serre fibrations as the relevant classes, cf. \cite{Hov}. Then $\Diag_\CC$ and $\WW(H,\Diag_\CC)$ are coloured $\Sg$-cofibrant (non-symmetric) operads, and Theorem \ref{main0} provides a model structure on the two categories $\Alg_\Top(\Diag_\CC)$ and $\Alg_\Top(\WW(H,\Diag_\CC))$ of $\CC$-diagrams of spaces, and of homotopy coherent such diagrams, respectively. In both cases, a map $f:X\to Y$ between diagrams is a weak equivalence if and only if $f(c):X(c)\to Y(c)$ is a weak equivalence in $\Top$ for every object $c$ in $\CC$. Since $\Top$ is left proper, Theorems \ref{comparison} and \ref{W}(i) imply that these categories are Quillen equivalent, under the Quillen adjunction$$\eps_!:\Alg_\Top(\WW(H,\Diag_\CC))\rlto\Alg_\Top(\Diag_\CC):\eps^*$$In particular, each homotopy coherent diagram $Y$ is weakly equivalent to a ``strict'' diagram $\eps^*(X)$; in other words, $Y$ can be rectified. This argument also applies if $\CC$ is an $\EE$-enriched category for which each hom-set $\CC(a,b)$ is a cofibrant object in $\EE$, and for which each unit map $I\to\CC(a,a)$ is a cofibration. (Here, $\EE$ is a cofibrantly generated monoidal model category with cofibrant unit and $H$ is an interval in $\EE$, all as before.) Write $\Diag_\EE(\CC)$ for the category of $\CC$-diagrams (i.e. algebras for the operad $\Diag_\CC$), and $\CohDiag_\EE(\CC)$ for the category of homotopy coherent $\CC$-diagrams (i.e. algebras for $\WW(H,\Diag_\CC)$). If $\EE$ satisfies the hypotheses of \ref{main0}, the categories $\Diag_\EE(\CC)$ and $\CohDiag_\EE(\CC)$ carry Quillen model structures (with weak equivalences and fibrations defined ``pointwise''). If moreover $\EE$ is left proper, then the adjunction$$\eps_!:\CohDiag_\EE(\CC)\rlto\Diag_\EE(\CC):\eps^*$$is a Quillen equivalence. In particular, the total left derived functor $L\eps_!$ provides a functorial rectification.

There is of course a result of the type above for every $\Sg$-cofibrant coloured operad $P$ (taking the place of $\Diag_\CC$). For easy reference, we state it explicitly.

\begin{cor}\label{rectification}(Rectification of homotopy $P$-algebras) Let $P$ be an admissible $\Sigma$-cofibrant operad in $\EE$ and assume that $\WW(H,P)$ is also admissible. Then, if moreover $\EE$ is left proper, the adjunction$$\eps_!:\Alg_\EE(P)\rlto\Alg_\EE(\WW(H,P)):\eps^*$$is a Quillen equivalence.\end{cor}

Another special case was pointed out to us by Ittay Weiss. Consider for a fixed set $\OO$ the operad $\Cat_\OO$ of Example \ref{openriched}. This is a $\Sg$-cofibrant non-symmetric coloured operad, and hence, if the hypotheses of \ref{main0} are met, the category of $\EE$-enriched categories with fixed set of objects carries a Quillen model structure, where a functor $F:A\to B$ between two such categories is a weak equivalence (resp. fibration) if and only if $F:A(x,y)\to B(x,y)$ is so for every $(x,y)\in\OO^2$. For an interval $H$, algebras for the operad $\WW(H,\Cat_\OO)$ are ``weak'' $\EE$-enriched categories, and it follows exactly as for diagrams that any such weak category is weakly equivalent to a strict one, if $\EE$ is left proper. 

For example, this applies if $\EE$ is the category of chain complexes with the projective model structure, where the algebras for $\WW(H,\Cat_\OO)$ are \emph{$A_\infty$-categories}, cf. Lyubashenko \cite{L}. It also applies if $\EE$ is the category $\Top$ of topological spaces where it yields the rectification result of Batanin \cite[Theorem 2.3]{Bat}. As another example, consider the category $\EE=\Cat$ of small categories with the ``folk'' model structure, for which the weak equivalences are the equivalences of categories and the cofibrations are one-to-one on objects. This is a left proper (cofibrantly generated) monoidal model category with the cartesian product as monoidal structure. Let $H$ be the groupoid on two objects generated by a single isomorphism between them. Then a $\Cat_\OO$-algebra is a $2$-category with $\OO$ as set of objects, while a $\WW(H,\Cat_\OO)$-algebra is an (unbiased) bicategory with the same objects. By Theorem \ref{comparison}, we obtain a Quillen equivalence of model categories \begin{center}($2$-categories on $\OO$)$\lra$(bicategories on $\OO$),\end{center} and in particular conclude (the known fact) that every bicategory is equivalent to a strict one.

For this special case $\EE=\Cat$, the rectification of diagrams (Corollary \ref{rectification}) specialises to the familiar fact that every ``pseudofunctor'' on a category $\CC$ is equivalent to a strict functor; or, equivalently, that every fibered category over $\CC$ is equivalent to a split one.

\section{Cosimplicial operads and weak morphisms of homotopy algebras}

Recall the fibered category of coloured operads for varying sets of colours, cf. \ref{changecol}. A \emph{cosimplicial operad} is by definition a cosimplicial object in this category. Such an object is given by a cosimplicial \emph{set} of colours $C^\bullet$, and for each $n$ a $C^n$-coloured operad $P^n$, together with maps$$P^n\lra\al^*(P^m)$$of $C^n$-coloured operads for any arrow $\al:[n]\to[m]$ in $\Delta$ with associated map (also denoted) $\al:C^n\to C^m$ between colour-sets. We will mostly denote such a cosimplicial operad by $P^\bullet$, leaving the colours implicit, and say that $P^\bullet$ is a cosimplicial operad \emph{over} $C^\bullet$. For a cosimplicial operad $P^\bullet$, the categories of algebras $\Alg_\EE(P^n),\,n\geq 0$, together form a (large) \emph{simplicial} category, which we denote by$$\Alg_\EE(P^\bullet).$$

If $P^\bullet$ is a cosimplicial operad over $C^\bullet$, then \emph{``geometric realisation''} with respect to $P^\bullet$ defines a functor$$\textrm{(simplicial sets)}\lra\textrm{(coloured operads)}$$sending a simplicial set $X$ to the operad $|X|_{P^\bullet}=X\otimes_\Delta P^\bullet$ with set of colours $|X|_{C^\bullet}=X\otimes_\Delta C^\bullet$. (Here we use that the fibered category of coloured operads is cocomplete.) An algebra $A$ over this operad $|X|_{P^\bullet}$ is given by an algebra $A_x$ over the $C^n$-coloured operad $P^n$, for each $n$-simplex $x\in X_n$, together with functorial isomorphisms$$\al^*(A_x)\overset{\cong}{\lra}A_{\al^*(x)}$$of $P^m$-algebras, where $\al:[m]\to[n]$ also denotes the induced map $P^m\to P^n$.

\subsection{Bimodules over $A_\infty$-monoids}Let $\Ass^0=\Ass$ be the operad whose algebras are associative unitary monoids. Let $\Ass^1=BiMod$ be the $3$-coloured operad whose algebras are triples $(A_0,M,A_1)$, where $A_0,A_1$ are $\Ass^0$-algebras, and $M$ is a left-$A_0$-right-$A_1$-bimodule, denoted\begin{diagram}[UO,small]A_0&\lTo^M&A_1.\end{diagram}$\Ass^0$ and $\Ass^1$ form part of a cosimplicial operad $\Ass^\bullet$ over the cosimplicial set $C^\bullet$ of colours, given by$$C^n=\{a_i\,|\,0\leq i\leq n\}\cup\{b_{ij}\,|\,0\leq i<j\leq n\},$$the colours for $n+1$ monoids $A_i$ and $\binom{n+1}{2}$ bimodules $M_{ij}$. The $\Ass^2$-algebras are diagrams of the form\begin{diagram}[UO,small,silent]A_0&\lTo^{M_{01}}&A_1\\&\luTo_{M_{02}}^{\overset{\theta}{\Longleftarrow}}&\uTo_{M_{12}}\\&&A_2\end{diagram}where $A_0,A_1,A_2$ are $\Ass^0$-algebras, $M_{ij}$ are bimodules, and $\theta:M_{01}\otimes M_{12}\to M_{02}$ is a map of $A_0$-$A_2$-bimodules, which is $A_1$-balanced in the sense that it equalises the two maps $M_{01}\otimes A_1\otimes M_{12}\dto M_{01}\otimes M_{12}$, (the tensor here refers to the monoidal structure of $\EE$). The $\Ass^3$-algebras are tetrahedra, whose faces are given by $\Ass^2$-algebras, and for which the induced square of bimodule maps commutes. The tensor product of two bimodules $M_{01}\otimes_{A_1}M_{12}$ is given as $\partial_1^*i_!$ for the maps of operads\begin{diagram}[UO,small]|\Lambda^{(1)}(2)|_{\Ass^\bullet}&\rTo^i&|\Delta(2)|_{\Ass^\bullet}&\lTo^{\partial_1}&|\Delta(1)|_{\Ass^\bullet}\end{diagram}and all the properties of this tensor product are actually consequences of properties of the Quillen pairs associated to such operads by geometric realisation of $\Ass^\bullet$.

Applying the Boardman-Vogt resolution, one gets a cosimplicial operad $\WW(\Ass^\bullet)$. The $\WW(\Ass^\bullet)$-algebras are $A_\infty$-algebras in $\EE$, the $\WW(\Ass^1)$-algebras are $\infty$-bimodules over such $A_\infty$-algebras, and so on. The map $\partial_1^*i_!$ for\begin{diagram}[UO,small]|\Lambda^{(1)}(2)|_{\WW(\Ass^\bullet)}&\rTo^i&|\Delta(2)|_{\WW(\Ass^\bullet)}&\lTo^{\partial_1}&|\Delta(1)|_{\WW(\Ass^\bullet)}\end{diagram}defines a tensor product $M\overset{\infty}{\otimes}N$ of such bimodules, and again, many properties of such a tensor product follow formally by considering Quillen pairs. For example, for $A_\infty$-algebras $A,B,C$ in $\EE$, and cofibrant $\infty$-bimodules $A\lTo^MB\lTo^NC$ and $A\lTo^{M'}B\lTo^{N'}C$, weak equivalences between $M$ and $M'$ and between $N$ and $N'$ induce a weak equivalence between $M\overset{\infty}{\otimes}_BN$ and $M'\overset{\infty}{\otimes}_BN'$. It seems worthwile to study $\infty$-bimodules in more detail from this point of view.

\subsection{Morphisms of $P$-algebras}Let $P$ be an ordinary (uncoloured) operad. Then $P=P^0$ is part of a cosimplicial operad $P^\bullet$, where $P^n$ is the operad whose algebras are $n$-simplices\begin{diagram}[UO,small]A_0&\rTo&A_1&\rTo&\cdots&\rTo&A_n\end{diagram}of morphisms of $P$-algebras. The set of colours of $P^n$ is $\{0,1,\dots,n\}$. For $n=1$, the operad $P^1$ has been described in Example \ref{morphisms}, and an explicit description of $P^n$ is very similar. In particular, if $P$ is $\Sg$-cofibrant, then so will be each $P^n$. The simplicial category $\Alg_\EE(P^\bullet)$ is the nerve of the category $\Alg_\EE(P)$ of $P$-algebras. As operad, $P^n$ can also be described as the pushout\begin{gather}\label{push}P^n=P^1\cup_{P^0}P^1\cup_{P^0}\cdots\cup_{P^0}P^1.\end{gather}More explicitly, for $n=2$ for example, $P^2$ is the pushout $P^1_{(01)}\cup_{P^0_{(1)}}P^1_{(12)}$ of operads coloured by $\{0,1,2\}$, where $P^1_{(01)}=\partial_{2!}(P^1)$ for $\partial_2:\{0,1\}\to\{0,1,2\},\,P^1_{(12)}=\partial_{0!}(P^1)$, and $P^0_{(1)}=\partial_{2!}\partial_{0!}(P)=\partial_{0!}\partial_{1!}(P)$ for $\partial_{2!}\partial_{0!}=\partial_{0!}\partial_{1!}:\{0\}\to\{0,1,2\}.$ In particular, $P^\bullet$ is ``$2$-coskeletal'' in the sense that for any simplicial set $X$, the inclusion $sk_2(X)\to X$ of the $2$-skeleton of $X$ induces an isomorphism\begin{gather}|sk_2(X)|_{P^\bullet}\overset{\cong}{\lra}|X|_{P^\bullet}\end{gather}of operads. In addition, the pushout (\ref{push}) for $n=2$ can be expressed by the isomorphism\begin{gather}|\Lambda^{(1)}[2]|_{P^\bullet}\overset{\cong}{\lra}|\Delta[2]|_{P^\bullet}.\end{gather}

\subsection{Weak morphisms of homotopy $P$-algebras}Let $P$ be a $\Sg$-cofibrant operad, and consider for each $n\geq 0$ the operad $P^n$ just introduced. Applying the Boardman-Vogt construction to each of them yields a cosimplicial operad$$\WW(P^\bullet)=\WW(H,P^\bullet).$$For $n=0$, the $\WW(P^0)=\WW(P)$-algebras are the homotopy $P$-algebras. For $n=1$, the $\WW(P^1)$-algebras are given by a pair $A_0,A_1$ of such homotopy $P$-algebras, together with a weak morphism between them, i.e. a morphism commuting with the operations of $P$ up to given (higher) homotopies. (For example, if $\EE$ is the category of differential graded vector spaces and $P=\Ass$, one essentially recovers the notion of an $A_\infty$-homomorphisms between $A_\infty$-algebras in this way.) The simplicial category $\Alg_\EE(\WW(P^\bullet))$ encodes in a sense the weak homomorphisms and all their compositions. To make this more precise, we shall prove that $\Alg_\EE(\WW(P^\bullet))$ looks like the nerve of a category up to Quillen equivalence, as expressed in the theorem below.

Consider for each $i=0,\dots,n-1$, the $n$ inclusions $\theta^i:[1]\to[n]$ mapping $0,1$ to $i,i+1$ respectively, and the induced maps $\theta^i:P^1\to P^n$ and $\WW(\theta^i):\WW(P^1)\to\WW(P^n)$. The identity (\ref{push}) above expresses $P^n$ as the pushout by the maps $\theta^i$, and from this it follows that the functor\begin{gather}\label{eq1}\Alg_\EE(P^n)\overset{\sim}{\lra}\Alg_\EE(P^1\cup_{P^0}\cdots\cup_{P^0}P^1)\end{gather}is an equivalence (even an isomorphism) of categories. The category of algebras over a pushout of operads is the pullback of the categories of algebras for the individual operads making up the pushout, so this equivalence (\ref{eq1}) can also be written as \begin{gather}\label{eq2}\Alg_\EE(P^n)\overset{\sim}{\lra}\Alg_\EE(P^1)\times_{\Alg_\EE(P^0)}\times\cdots\times_{\Alg_\EE(P^0)}\Alg_\EE(P^1).\end{gather}The maps $\WW(\theta^i):\WW(P^1)\to\WW(P^n)$ together induce a map\begin{gather}\theta:\WW(P^1)\cup_{\WW(P^0)}\cdots\cup_{\WW(P^0)}\WW(P^1)\lra\WW(P^n)\end{gather}which induces an equivalence of homotopy categories. More precisely,
\begin{thm}Let $\EE$ be a left proper (cofibrantly generated) monoidal model category and $P$ a $\Sg$-cofibrant operad in $\EE$ for which the associated operads $\WW(P^n)$ all are admissible. Then for each $n\geq 2$, the map $\theta$ induces a Quillen equivalence\begin{gather}\label{eq3}\theta_!:\Alg_\EE(\WW(P^1)\cup_{\WW(P^0)}\cdots\cup_{\WW(P^0)}\WW(P^1))\rlto\Alg_\EE(\WW(P^n)):\theta^*,\end{gather}or, equivalently, a Quillen equivalence of model categories \begin{gather}\label{eq4}\Alg_\EE(\WW(P^1))\times_{\Alg_\EE(\WW(P^0))}\times\cdots\times_{\Alg_\EE(\WW(P^0))}\Alg_\EE(\WW(P^1))\overset{\sim}{\leftrightarrow}\Alg_\EE(\WW(P^n)).\end{gather}\end{thm}

\begin{rmk}Before we embark on the proof, it is good to be more explicit about pullbacks of categories like the one in (\ref{eq4}). For functors $\phi:\BB\to\AA\leftarrow\CC:\psi$, one can consider the ``strict'' pullback category $\BB\times_\AA\CC$ where objects are pairs $(B,C)$ with $\phi (B)=\psi (C)$, and one can consider the ``pseudo-''pullback whose objects are triples $(B,C,\xi)$ where $\xi:\phi( B)\overset{\cong}{\lra}\psi( C)$ is an isomorphism in $\AA$. These constructions are in general different, and the latter construction is often more useful. However, if one of the functors, say $\phi$, is a ``categorical fibration'' in the sense that any isomorphism $\phi(B)\overset{\cong}{\lra}A$ in $\AA$ is the image under $\phi$ of an isomorphism $B\overset{\cong}{\lra}B'$ in $\BB$, then the strict pullback category and the pseudo-one are equivalent, and there is no need to distinguish between the two. This is so in the particular case at hand. In fact, for any map $Q\to R$ between operads, the induced functor $\Alg_\EE(R)\to\Alg_\EE(Q)$ is a categorical fibration.\end{rmk}

\emph{Proof of Theorem.} We treat the case $n=2$ only. (The argument is the same for higher $n$, but the notation is more involved.) Consider the cube

\begin{gather}\label{cube1}\begin{diagram}[UO,small,silent]P^1\cup_{P^0}P^1&&\lTo&&P^1&&\\&\luTo&&&\uTo&\luTo&\\\uTo&&\WW(P^1)\cup_{\WW(P^0)}\WW(P^1)&\lTo&\HonV&&\WW(P^1)\\&&\uTo&&\vLine&&\\P^1&\lTo&\VonH&\hLine&P^0&&\uTo\\&\luTo&&&&\luTo&\\&&\WW(P^1)&&\lTo&&\WW(P^0)\end{diagram}\end{gather}and the associated diagram of categories of algebras. The latter diagram is of the following type
\begin{gather}\label{cube2}\begin{diagram}[UO,small,silent]\BB\times_\AA\CC&&\rTo^{q^*}&&\CC&&\\&\rdTo&&&\vLine&\rdTo&\\\dTo^{p^*}&&\BB'\times_{\AA'}\CC'&\rTo^{q'^*}&\HonV&&\CC'\\&&\dTo_{p'^*}&&\dTo^{r^*}&&\\\BB&\hLine&\VonH&\rTo^{u^*}&\AA&&\dTo_{r'^*}\\&\rdTo_{\be^*}&&&&\rdTo^{\al^*}&\\&&\BB'&&\rTo_{u'^*}&&\AA'\end{diagram}\end{gather}Here $\AA,\BB,\CC$ and $\AA',\BB',\CC'$ are Quillen model categories, and the (pseudo)pullback categories $\BB\times_\AA\CC,\BB'\times_{\AA'}\CC'$ also carry a model structure, for which the projection functors ($p^*,q^*$ and $p'^*,q'^*$) preserve fibrations and weak equivalences. In fact, all functors are induced by operad maps and hence are right parts of Quillen pairs. Furthermore, by Theorem \ref{comparison}, the pairs $(\al_!,\al^*),(\be_!,\be^*)$ and $(\gamma_!,\gamma^*)$ are Quillen equivalences. We now need to show that the induced functor$$(\be^*,\gamma^*):\BB\times_\AA\CC\lra\BB'\times_{\AA'}\CC'$$is also part of a Quillen equivalence. For this, it is enough to show that the left adjoints $\be_!$ and $\gamma_!$ together define a functor $(\be_!,\gamma_!)$ from $\BB\times_\AA\CC$ into $\BB'\times_{\AA'}\CC'$. Because then this functor will be left adjoint to $(\be^*,\gamma^*)$, and the appropriate units and counits will be weak equivalences. In other words, we need to prove that for objects $B'$ of $\BB'$ and $C'$ of $\CC'$, and a given isomorphism $u'^*B'\cong v'^*C'$, there is an induced natural isomorphism $v^*\be_!B\cong u^*\gamma_!C$. This would indeed follow if the two squares on the bottom and the right of the cube (\ref{cube2}) satisfy the ``Beck-Chevalley condition'' (or ``projection formula''), stating that the natural maps$$\al_!u'^*\lra u^*\be_!\text{ and }\al_!v'^*\lra v^*\gamma_!$$are isomorphisms. In the particular case at hand for the two squares of operads $(i=0,1)$\begin{diagram}[UO,small,silent]P^1&\lTo^{\partial^i}&P^0\\\uTo_{\eps^1}&&\uTo^{\eps^0}\\\WW(P^1)&\lTo_{\WW(\partial^i)}&\WW(P^0)\end{diagram}these natural maps are the maps $(i=0,1)$ $$\eps^0_!\WW(\partial^i)^*\lra(\partial^i)^*\eps^1_!$$on the categories of algebras. To see that they are isomorphisms, we note that for a $P^1$-algebra $A=(A_0\overset{f}{\lra}A_1)$, one has$$(\partial^0)^*A=A_1\text{ and }(\partial^1)^*(A)=A_0$$so that for a $P^0$-algebra $B$,$$(\partial^0)_!(B)=(0\lra B)\text{ and }(\partial^1)_!(B)=(B\overset{id}{\lra}B).$$Moreover, the functors $(\partial^0)^*$ and $(\partial^1)^*$ also have right adjoints, given by$$(\partial^0)_*(B)=(B\overset{id}{\lra}B)\text{ and }(\partial^1)_*(B)=(B\lra T)$$where $T$ is the terminal object of $\EE$ (and hence of $\Alg_\EE(P)$). A similar description applies to the functors induced by $\WW(\partial^i)$. For example, for $i=0$, given a $\WW(P^0)$-algebra $B$, the two objects $B_0$ and $B_1$ given by $B_0=B=B_1$ together form a $\WW(P^1)$-algebra in a canonical way. (A tree representing an operation in $\WW(P^1)$ acts, by first replacing all the colours ($0$ or $1$) on the edges of the tree by $1$, and then letting the tree act via the given $\WW(P^0)$-algebra structure on $B$.) This gives the right adjoint $\WW(\partial^0)_*(B)$. Similarly, for $i=1$, and a given $\WW(P^0)$-algebra $B$, the two objects $B_0=B$ and $B_1=T$ (terminal object) form a $\WW(P^1)$-algebra, and this is the value of the right adjoint $\WW(\partial^1)_*(B)$. From this explicit description, we see that there is a natural isomorphism$$\WW(\partial^i)_*(\eps^0)^*(B)\cong(\eps^1)^*(\partial^i)_*(B)$$for any $P^0$-algebra $B$. By taking left adjoints of the functors involved, it follows that there is a natural isomorphism$$\eps^0_!\WW(\partial^i)^*(A)\cong(\partial^i)^*\eps^1_!(A)$$for any $\WW(P^1)$-algebra $A$, as was to be shown.\qed

\begin{rmk}Consider again the cube (\ref{cube1}) in the proof above. The pushouts of operads are calculated by first pushing the operads forward to the set of colours $\{0,1,2\}$, and then taking the pushout of $\{0,1,2\}$-coloured operads. In other words, we could have written a diagram of $\{0,1,2\}$-coloured operads, whose bottom face, for example, is\begin{diagram}[UO,small]\partial_{2!}(P^1)&\lTo&\partial_{2!}\partial_{0!}(P^0)\\\uTo&&\uTo\\\partial_{2!}(\WW(P^1))&\lTo&\partial_{2!}\partial_{0!}(\WW(P^0)).\end{diagram}Since $\partial_{0!}P^0\to P^1$ is a $\Sg$-cofibration between $\Sg$-cofibrant $\{0,1\}$-coloured operads, and since $\partial_{i!}$ commutes with the Boardman-Vogt construction, cf. Remark \ref{changecol2}, it follows from Theorem \ref{W}(ii) that the map \begin{diagram}[UO,small]\partial_{2!}(\WW(P^1))&\lTo&\partial_{2!}\partial_{0!}(\WW(P^0))\end{diagram} is a cofibration of $\{0,1,2\}$-coloured operads. In other words, deleting the $\partial_{i!}$ from the notation but interpreting the diagram (\ref{cube1}) as a diagram of $\{0,1,2\}$-coloured operads, the right hand front vertical map is a cofibration, and hence so is the left hand front vertical map (by pushout). Now, inscribe in the right and left hand faces of the cube the pushouts along these two cofibrations of the right and left bottom weak equivalences $\WW(P^0)\to P^0$ and $\WW(P^1)\to P^1$, respectively. \emph{If} the model structure on $\{0,1,2\}$-coloured operads is left proper, the pushouts are again weak equivalences. It then follows by the ``2-out-of-3'' axiom for weak equivalences and Lemma 6.9 of \cite{BM2} that the map $\WW(P^1\cup_{\WW(P^0)}\WW(P^1)\lra P^1\cup_{P^0}P^1$ is also a weak equivalence. Hence, by commutativity of the square\begin{diagram}[UO,small]\WW(P^1)\cup_{\WW(P^0)}\WW(P^1)&\rTo&P^1\cup_{P^0}P^1\\\dTo&&\dTo^{\cong}\\\WW(P^2)&\rTo&P^2,\end{diagram}we find that $\WW(P^1)\cup_{\WW(P^0)}\WW(P^1)\to\WW(P^2)$ is a weak equivalence as well. In this way, we can improve on the previous theorem:\end{rmk}
\begin{thm}\label{Segal}Under the hypotheses of the previous theorem, assume in addition that the model structure on coloured operads is left proper. Then for each $n$, the map$$\WW(P^n)\lra\WW(P^1)\cup_{\WW(P^0)}\cdots\cup_{\WW(P^0)}\WW(P^1)$$is a weak equivalence of operads.\end{thm}
\begin{rmk}Under the circumstances of the previous theorem, one can construct a homotopy category of homotopy $P$-algebras and homotopy classes of weak maps, denoted $\Ho(\WW(P^\bullet))$, as follows. 

The objects are the cofibrant-fibrant $\WW(P^0)$-algebras $A=(A,\al)$, where $\al:\WW(P^0)\to\End(A)$ denotes the structure map. Given two such objects $(A,\al)$ and $(B,\be)$, let $\End(A,B)$ denote the $\{0,1\}$-coloured endomorphism-operad with $A$ concentrated in colour $0$, and $B$ in colour $1$. Then an arrow $F:(A,\al)\to(B,\be)$ in $\Ho(\WW(P^\bullet))$ is a map $\WW(P^1)\to\End(A,B)$ in $\Ho((\WW(P^0)+\WW(P^0))/\Oper_{\{0,1\}}(\EE))$, the homotopy category associated to the model category of $\{0,1\}$-coloured operads under $\WW(P^0)+\WW(P^0)$, where the first summand is concentrated in colour $0$ and the second in colour $1$. 

Now, $\WW(P^0)+\WW(P^0)\to\End(A,B)$ is a fibrant object in this model category, and $\WW(P^0)+\WW(P^0)\to\WW(P^1)$ is cofibrant (by Theorem \ref{W}(ii)), so a map $F:(A,\al)\to(B,\be)$ in $\Ho(\WW(P^\bullet))$ is represented by a map $f$ making the diagram\begin{diagram}[UO,small,silent]\WW(P^0)+\WW(P^0)&\rTo^{\al+\be}&\End(A,B)\\\dTo&\ruTo_f&\\\WW(P^1)\end{diagram}commute, and two such maps $f$ and $f'$ represent the same arrow if and only if they are homotopic rel. $\WW(P^0)+\WW(P^0)$.

To compose two such maps $F:(A,\al)\to(B,\be)$ and $G:(B,\be)\to(C,\gamma)$ represented by $f$ and $g$ respectively, we consider the model category of $\{0,1,2\}$-coloured operads under $|sk_0(\Delta(2))|_{\WW(P^\bullet)}=\WW(P^0)+\WW(P^0)+\WW(P^0)$. In this category, the map $|\Lambda^{(1)}(2)|_{\WW(P^\bullet)}\lra|\Delta(2)|_{\WW(P^\bullet)}$, i.e.$$\WW(P^1)\cup_{\WW(P^0)}\WW(P^1)\lra\WW(P^2),$$is a weak equivalence by Theorem \ref{Segal}, and the objects involved are cofibrant. Indeed, the map $\WW(P^0)+\WW(P^0)+\WW(P^0)\to\WW(P^2)$ is a cofibration by \ref{W}(ii), and $\WW(P^0)+\WW(P^0)+\WW(P^0)\to\WW(P^1)\cup_{\WW(P^0)}\WW(P^1)$ is a cofibration because it is a pushout of the sum of two copies of the cofibration $\WW(P^0)+\WW(P^0)\to\WW(P^1)$. 

Thus $f$ and $g$ together define a map $\WW(P^1)\cup_{\WW(P^0)}\WW(P^1)\to\End(A,B,C)$ under $\WW(P^0)+\WW(P^0)+\WW(P^0)$, and this map extends uniquely up to homotopy (rel. $\WW(P^0)+\WW(P^0)+\WW(P^0)$) to a map $h:\WW(P^2)\to\End(A,B,C)$. The composition of $h$ with $\partial^1:\WW(P^1)\to\WW(P^2)$ gives a map $\WW(P^1)\to\End(A,C)$ of $\{0,1\}$-coloured operads under $\WW(P^0)+\WW(P^0)$, representing the composition $GF:(A,\al)\to(C,\gamma)$ in $\Ho(\WW(P^\bullet))$. It follows by purely model theoretic facts (similar to the one used to define $h$) that this composition is well defined on homotopy classes, and is associative.\end{rmk}

\section{Appendix: Coloured operads as monoids}

In this appendix, we will explicitly describe a tensor product on coloured collections for which coloured operads are monoids. This tensor product specialises to the well-known (Smirnov) tensor product on ordinary collections if there is just one colour. We continue to work with a cocomplete closed symmetric monoidal category $\EE=(\EE, \otimes,I)$.

\subsection{The categories $\mathbb{F}(C), \mathbb{F}^\circ(C),\mathbb{F}^\leq(C)$}Let $C$ be a fixed set of colours. We consider the category $\mathbb{F}(C)$ whose objects are triples $(S,s_0,\al)$ where $S$ is a finite set, $s_0{\in S}$ is a base point, and $\al:S\rightarrow C$ is a function.  We will often simply denote such an object by $S$. An arrow $(S,s_0,\al)\stackrel{\sg}{\rightarrow}(T,t_0,\beta)$ in this category is a basepoint preserving bijection $\sg:S \stackrel{\sim}{\rightarrow} T$ for which $\beta \sg = \al$.

Each finite pointed set is isomorphic to a set of the form $\{1,\dots,n,\star\}$ with $n\geq 0$ and $\star \notin\{1,\dots,n\}$ viewed as the basepoint.  Let $\mathbb{F}^\circ(C)$ be the full subcategory on $\mathbb{F}(C)$ given by those finite sets.  The category $\mathbb{F}^\circ(C)$ is a sum of translation groupoids$$(C^n\times C)\rtimes \Sg_n\hspace*{1in}(n\geq 0)$$where $\Sg_n$ acts on $C^n\times C$ by permuting the first $n$ coordinates, $(c_1,\dots,c_n;c)^\sg=(c_{\sg(1)},\dots,c_{\sg(n)};c)$.  The inclusion $\mathbb{F}^\circ(C) \rightarrow \mathbb{F}(C)$ is an equivalence of categories. Let us suppose that $C$ is equipped with a \emph{linear order} $\leq$.  Then there is a smaller full subcategory $\mathbb{F}^\leq(C)\subseteq \mathbb{F}^\circ(C)$, still equivalent to $\mathbb{F}(C)$, given by only those objects $\al:\{1,\dots,n,\star\} \rightarrow C$ for which $\al(1)\leq \dots \leq \al(n)$.  The category $\mathbb{F}^\leq(C)$ is a sum of symmetric groups:  Indeed, for each sequence $c_1 \leq \dots \leq c_n$ in $C$, write $\Sg_{c_1\dots c_n}\subseteq \Sg_n$ for the subgroup of permutations $\sg$ for which $c_{\sg(1)}\leq \dots \leq c_{\sg(n)}$.  (Note that this implies that $(c_1,\dots,c_n)=(c_{\sg(1)},\dots,c_{\sg(n)})$.  Then$$\mathbb{F}^\leq(C)=\coprod_c \coprod_{c_1\leq \dots \leq c_n}\textstyle \sum_{c_1\dots c_n}.$$Since any set $C$ can always be given a linear order, the category $\mathbb{F}(C)$ is equivalent to one of the form $\mathbb{F}^\leq(C)$.

\subsection{$C$-collections}A $C$-collection is a functor$$K:\mathbb{F}(C)^{op}\rightarrow \EE$$and a map of $C$-collections is a natural transformation.  This defines a category$$Coll_C(\EE)$$of $C$-coloured collections in $\EE$.  The monoidal structure of $\EE$ induces a pointwise monoidal structure on $Coll_C(\EE)$,$$(K \otimes L)(S)=K(S)\otimes L(S)$$where $S=(S,s_0,\al)$ is any object in $\mathbb{F}(C)$.  There is another tensor product on $Coll_C(\EE)$, parametrising the composition of coloured operads, and defined as follows.  Consider two objects $(S,s_0,\al)$ and $(R,r_0,\gamma)$ in $\mathbb{F}(C)$.  For a map $(R,r_0)\rightarrow (S,s_0)$ of finite pointed sets (not necessarily a bijection), write$$R_s=f^{-1}(s)\cup\{s\}$$where we assume $s \notin f^{-1}(s)$, and view $s$ as the base point of $R_s$.  Write $\gamma_s:R_s\rightarrow C$ for the map defined by restricting $\gamma$ to $f^{-1}(s)$ and by $\gamma_s(s)=\al(s)$.  Now define
\begin{gather}\label{colimit}(K \square L)(R,r_0,\gamma)=\cli\,[K(S,s_0,\al)\otimes \bigotimes_{s \in S-\{s_0\}}L(R_s,s,\gamma_s)],\end{gather}where the colimit is taken over the category whose objects are factorisations of $\gamma$ as a composition $\gamma=\al f$,and whose arrows are commutative diagrams
\begin{gather}\label{arrow}\begin{diagram}[UO,small]&&(S,s_0)&&\\&\ruTo^f&&\rdTo^\al&\\(R,r_0)&&\dTo^\sg&&C\\&\rdTo_g&&\ruTo_\beta&\\&&(T,t_0)&&\end{diagram}\end{gather}where $\sg$ is an arrow in $\mathbb{F}(C)$, i.e. a bijection.  The expression ``$\displaystyle{\bigotimes_{s \in S-\{ s_0\}}}$'' denotes the set-indexed tensor product; it can be defined for any finite family $\{E_j\}_{j \in J}$ as$$\bigotimes_{j \in J}E_j=\left(\coprod_{\sg:\{1,\dots,n\}\rightarrow J}(E_{\sg(1)}\otimes \dots \otimes E_{\sg(n)})\right)_{\Sg_n}$$where $\Sg_n$ acts by composition on bijections $\sg$ and by canonical isomorphisms of $\EE$.  Any arrow $(\ref{arrow})$ then induces a map$$K(T,t_0, \beta)\otimes\bigotimes_{t\neq t_0}L(R_t,t,\gamma_t)\stackrel{\sg^\star}{\rightarrow}K(S,s_0,\al)\otimes \bigotimes_{s \neq s_0}L(R_s,s,\gamma_s)$$and the colimit (\ref{colimit}) is along these maps. This tensor product is associative but not commutative (and only closed on one side). It has a 2-sided unit $U$, given by$$U(S,s_0,\al)=\begin{cases}I\text{ whenever}&|S|=2\,\text{ and }\al\text{ is constant};\\ 0&\text{otherwise}.\end{cases}$$

The category of collections $Coll_C(\EE)$ is of course equivalent to each of the smaller categories$$Coll^\circ_C(\EE)\text{ and }Coll^\leq_C(\EE)$$of contravariant functors on $\mathbb{F}^\circ(C)$ respectively on $\mathbb{F}^\leq(C)$.  An object $X$ of $Coll^\circ_C(\EE)$ assigns to each sequence $(c_1,\dots,c_n;c)$ an object $X(c_1,\dots,c_n;c)$ of $\EE$ and to each $\sg \in \Sg_n$ a map $\sg^\star:X(c_1,\dots,c_n;c)\rightarrow X(c_{\sg(1)},\dots,c_{\sg(n)};c)$, functorial in $\sg$.  The tensor product $\square$ takes the familiar form$$(X \square Y)(d_1,\dots,d_k;c)=$$
$$\coprod_n\left[\coprod_{c_1\dots c_n}\coprod_{k_1 + \dots +k_n=k}X(c_1,\dots,c_n;c)\otimes\left(\bigotimes^n_{i=1}Y(d^{(i)};c_i)\right)\right]\otimes_{\Sg_n \ltimes (\Sg_{{k}_1}\times \dots \times \Sg_{{k}_n})}I[\Sg_n]$$where $d^{(i)}=(d_{{k}_{1+\dots +}k_{i-1^{+1}}},\dots,d_{{k}_1+\dots +k_i})$.  An object $X$ of $Coll^\leq_C(\EE)$ only assigns an object $X(c_1,\dots, c_n;c)$ to a sequence for which $c_1 \leq \dots \leq c_n$, and assigns a map $\sg^\star$ only to $\sg \in \Sg_{{c}_1 \dots c_n}$.  In other words$$Coll^\leq_C(\EE)=\prod_{{c}_1 \leq \dots \leq c_n}\EE^{\Sg_{{c}_1 \dots c_n}}$$(This resembles most the definition of uncoloured collections as $Coll(\EE)=\prod_n \EE^{\Sg_n})$.  On $Coll^\leq_C(\EE)$, the tensor product $\square$ takes the form$$(X \square Y)(d_1,\dots, d_k;c)=\displaystyle\coprod_n[\dots]_{{\Sg}_{{c}_1\dots c_n}}$$where $[\dots]$  is$$\coprod_{{c}_1 \leq \dots \leq c_n}\coprod_{{k}_1+ \dots + k_n=k}\coprod_\rho X(c_1,\dots,c_n;c)\otimes\bigotimes^n_{i=1}Y(d_\rho^{(i)}; c_i)$$and $\rho$ ranges over all permutations for which $d \cdot \rho$ is the concatenation of $n$ ordered sequences $d_\rho^{(1)},\dots,d_\rho^{(n)},$ where $d=(d_1,\dots,d_k)$. The group $\Sigma_k$ acts as follows: if $d \cdot \tau = d$ then $\tau$ acts by mapping the summand for $\rho$ to that for $\tau^{-1}\rho \tau,$ and by acting on each of the blocks $d_\rho^{(i)}$ by $\rho \tau\rho^{-1}$ (this makes sense because $d_{\tau (\ell)}= d_\ell$ for all $\ell$).

 \subsection{Operads}For a set of colours $C$, one can now define the category of $C$-coloured operads as the category of (unitary and associative) monoids in $Coll_C(\EE)$ equipped with the $\square$-product.  By the equivalences of (monoidal) categories $$Coll_C(\EE) \simeq Coll^\circ_C(\EE) \simeq Coll^\leq_C(\EE),$$ an equivalent definition of $C$-coloured operads is as monoids in each of the other two categories.  The definition we presented in$\S$ 1 is that of a monoid in $Coll^\circ_C(\EE)$.  The most convenient definition, however, is that of a monoid in$Coll_C^\leq(\EE)$, cf. Remark \ref{order}. Such a monoid $P$ is given by an object $P(c_1,\dots,c_n;c)$ for each ordered sequence $c_1 \leq \dots\leq c_n$ and each $c$, a right action by $\Sg_{{c}_1 \dots c_n}$ on that object, unit $I \rightarrow P(c;c)$ for each $c$, and composition maps. To describe these, let us introduce some notation:  The group $\Sg_k$ acts on sequences $d=(d_1,\dots,d_k)$ on the right, by $(d \cdot\sg)_i=d_{\sg(i)} $.  Also, if $d^{(1)},\dots,d^{(n)}$ are $n$ sequences, say $d^{(i)}= (d^{(i)}_1,\dots,d^{(i)}_{{k}_i})$ of length $k_i$, we write $d=(d^{(1)},\dots,d^{(n)})$ for the concatenated sequence of length $k=k_1+\dots +k_n$.  (Each of the $d^{(i)}$ can be ordered, but this does not imply that $d$ is, of course.) Now, for each $c$, and each $n$-tuple of ordered sequences $d^{(1)},\dots,d^{(n)}$, and each $\rho\in \Sg_k$ for which $d \cdot \rho$ is ordered, $\rho$ is equipped with a composition map$$P(c_1,\dots ,c_n;c)\otimes \bigotimes^n_{i=1}P(d^{(i)};c_i)\stackrel{\gamma_\rho}{\rightarrow}P(d \cdot \rho;c),$$and besides being unitary and associative, these should satisfy two equivariance conditions:

(i)  For each $\sg \in \Sg_{{c}_1\dots c_n}$, the diagram\begin{diagram}P(c_1, \dots, c_n;c)\otimes \bigotimes^n_{i=1}P(d^{(i)};c_i)&\rTo^{\gamma_\rho}&P(d\cdot \rho;c)\\\dTo^{\sg^\star \otimes \tilde{\sg}}&& \dTo_{(\rho^{-1}\bar{\sg}\tau)^\star}\\P(c_1,\dots,c_n;c) \otimes \bigotimes^n_{i=1}P(d^{\sg(i)};c_{\gamma(i)})&\rTo^{\gamma \tau}&P(d \cdot \bar{\sg}\cdot \tau;c)\end{diagram} commutes. Here $\bar{\sg}\epsilon \Sg_k$ denotes the block permutation given by $\sg \in \Sg_n$ and $k=k_1+ \dots + k_n$, and $\tilde{\sg}$ is the canonical isomorphism in $\EE$; furthermore, $\rho$ and $\tau$ are any permutations putting $d=(d^{(1)},\dots,d^{(n)})$ and $d\cdot\bar{\sg}=(d^{\sg(1)},\dots,d^{\sg(n)})$ in the right order, so that $d \cdot \rho$ and $(d \cdot \bar{\sg})\cdot\tau$ are the same sequence, and $\rho^{-1}\bar{\sg}\tau \in \Sg_{d\cdot\rho}$.

Note that this condition is non-vacuous when $\sg$ is the identity, where it says that for $\rho$ as in (5) and $\xi \in\Sg_{d \cdot \rho}$,$$\xi^\star \circ \gamma_\rho = \gamma_{\rho \xi}$$

(ii)  For each $n$-tuple of permutations $\sg_i \in\Sg_{{d}^{(i)}}$ and each $\rho \in \Sg_d$, the diagram\begin{diagram}P(c_1,\dots,c_n;c) \otimes \bigotimes^n_{i=1}P(d^{(i)};c_i) &\rTo^{\gamma_\rho}&P(d \cdot \rho; c)\\\dTo^{id \otimes \sg_1^\star \otimes \dots \otimes \sg^\star_n}&& \dTo_{(\rho^{-1} \sg \rho)^\star}\\P(c_1,\dots,c_n;c) \otimes \bigotimes^n_{i=1} P(d^{(i)};c)&\rTo^{\gamma_\rho}&P(d\cdot \rho;c)\end{diagram} commutes, where $\sg=\sg_1 \times \dots \times \sg_n \in \Sg_{{d}^{(1)}} \times \dots\times \Sg_{{d}^{(n)}}$.  (We could also write this as a commutative upper triangle, with dotted arrow $\gamma_{\sg^{-1} \cdot \rho})$. 

Thus, we obtain three equivalent categories of $C$-coloured operads,$$Oper_C(\EE)\simeq Oper^\circ_C(\EE) \simeq Oper^\leq_C(\EE),$$as monoids in $Coll_C(\EE), Coll_C^\circ (\EE)$ and $Coll^\leq_C(\EE)$, respectively.

\vspace{5ex}

\noindent{\sc Universit\'e de Nice, Laboratoire J.-A. Dieudonn\'e, Parc Valrose, 06108 Nice Cedex, France.}\hspace{2em}\emph{E-mail:} cberger$@$math.unice.fr\vspace{2ex}

\noindent{\sc Mathematisch Instituut, Universiteit Utrecht, The Netherlands,} and {\sc Department of Pure Mathematics, University of Sheffield, UK.}\\\emph{E-mail:} moerdijk$@$math.uu.nl

\end{document}